\documentclass[preprint,12pt]{elsarticle}
\usepackage[a4paper,margin=27mm]{geometry}
\usepackage{amsmath,amssymb,amsthm,mathtools}
\usepackage{xcolor}
\usepackage[hidelinks]{hyperref}
\newtheorem{theorem}{Theorem}[section]
\newtheorem{lemma}[theorem]{Lemma}
\newtheorem{corollary}[theorem]{Corollary}
\newtheorem{proposition}[theorem]{Proposition}
\theoremstyle{definition}
\newtheorem{assumption}[theorem]{Assumption}
\theoremstyle{remark}
\newtheorem{remark}[theorem]{Remark}
\newcommand{\PP}{\mathbb P}
\newcommand{\NN}{\mathbb N_0}
\newcommand{\supp}{\operatorname{supp}}
\newcommand{\diam}{\operatorname{diam}}
\newcommand{\Span}{\operatorname{span}}

\usepackage{graphicx,enumitem,booktabs,tabularx,flafter}
\theoremstyle{definition}

\DeclareMathOperator{\ran}{range}

\newcommand{\norma}[1]{\lVert#1\rVert_a}

\begin{document}
	\begin{frontmatter}
		\title{Product Enrichment of Spline Spaces for IGA:\protect\\Polynomial Reproduction, Approximation, and Stable Computation}

	\author{Shaojie Luo, Yuanpeng Zhu*} 
\cortext[cor]{Corresponding author. \newline
	E-mail addresses: 202510188359@mail.scut.edu.cn (S.J. Luo), ypzhu@scut.edu.cn (Y.P. Zhu)}
	\address{School of Mathematics and Statistics, Central South University, Changsha, $410083$, PR China}

\begin{abstract}
We study product enrichment of spline spaces under a fixed isoparametric
geometry map for isogeometric analysis (IGA). For each fixed integer $\ell\ge1$, a locally
$L^\infty$-stable initial operator reproducing $\PP_k$ yields an explicit
linear map from $L^2(D)$ into $W_h^\ell$ that reproduces $\PP_{\ell k}$
and is locally stable in $L^2$ and $L^\infty$, uniformly in the mesh scale.
Under additional initial-space regularity and derivative stability,
the construction gives integer-order Sobolev error estimates of order
$h^{s-m}$ for $1\le s\le\ell k+1$ and admissible $m$.
Independently, the geometry representation implies reproduction of
physical polynomials of degree at most $\ell$ by the enriched space.
Polynomial containment also yields algebraic and analytic best-approximation
estimates as the product order increases. For every fixed product order,
we give a compatible redundant Galerkin formulation and an explicit
initial-space injection. A local spectral construction and an initial-space
correction give a two-level energy condition-number bound without requiring
independent product generators. These solver results concern exact
arithmetic; coarse-space cost, quadrature consistency, boundary constraints,
and numerical rank truncation require separate control.
\end{abstract}

\begin{keyword}
	{spline approximation;  product enrichment; quasi-interpolation; polynomial reproduction; isogeometric analysis}

	
\end{keyword}
\end{frontmatter}

\section{Introduction}\label{sec:intro}
Isogeometric analysis represents the geometry and the discrete field using
related approximation spaces. Its approximation theory therefore depends
both on the parameter-space construction and on the regularity of the
geometry map~\cite{Bazilevs2006}. Polynomial reproduction in parameter
coordinates must be distinguished from polynomial reproduction in physical
coordinates. This distinction motivates the product enrichment studied here.

Let $d\ge1$, $D=(0,1)^d$, and $\NN=\{0,1,2,\ldots\}$. For a multi-index
$\gamma\in\NN^d$, set $|\gamma|=\sum_{\nu=1}^d\gamma_\nu$ and
$z^\gamma=\prod_{\nu=1}^d z_\nu^{\gamma_\nu}$. We write
\[
\PP_q=\Span\{z^\gamma:|\gamma|\le q\},\qquad
\PP_q(E)=\{p|_E:p\in\PP_q\}.
\]
Thus $\PP_q$ denotes complete polynomials of total degree at most $q$.
Unless stated otherwise, all function identities are understood almost
everywhere. For integer $j\ge1$, $|u|_{H^j(E)}^2$ is the sum of the
squared $L^2(E)$ norms of weak derivatives of total order $j$; we set
$|u|_{H^0(E)}=\|u\|_{L^2(E)}$ and use the corresponding full Sobolev norm.

For $h\in(0,1]$, let $\mathcal I_h=\{1,\ldots,N_h\}$ and let
$\{B_{i,h}\}_{i\in\mathcal I_h}\subset L^\infty(D)$ be linearly
independent real-valued functions satisfying
\begin{equation}\label{eq:partition_unity}
\sum_{i\in\mathcal I_h}B_{i,h}=1.
\end{equation}
Set
\begin{equation}\label{eq:initial-space}
W_h^1=W_h=\Span\{B_{i,h}:i\in\mathcal I_h\},\qquad \dim W_h^1=N_h.
\end{equation}
The abstract approximation construction does not require polynomial local
pieces. Polynomial spline structure will be imposed separately for the
computational realization.

Let $\Omega\subset\mathbb R^d$ be a bounded domain. A fixed geometry
map is represented, at every mesh level under consideration, by
\begin{equation}\label{eq:ext_geo}
F:D\longrightarrow\Omega,\qquad
F(\xi)=\sum_{i\in\mathcal I_h}B_{i,h}(\xi)X_{i,h},
\qquad X_{i,h}\in\mathbb R^d.
\end{equation}
For algebraic reproduction, $F$ is taken to be a continuous bijection;
the Sobolev results impose bi-Lipschitz continuity and the stated higher
regularity. With $F^*v=v\circ F$, the initial physical space is
\[
V_h^1=\{w_h\circ F^{-1}:w_h\in W_h^1\}.
\]
The partition of unity and the geometry representation imply
$F^*\PP_1(\Omega)\subseteq W_h^1$. However, for coordinate indices
$1\le\mu,\nu\le d$,
\begin{equation}\label{eq:proper_sub}
(x_\mu x_\nu)\circ F
=F_\mu F_\nu
=\sum_{i,j\in\mathcal I_h}X_{i,h}^{\mu}X_{j,h}^{\nu}B_{i,h}B_{j,h},
\end{equation}
and these products need not belong to $W_h^1$. For example, take
$D=\Omega=(0,1)$, $W_h^1=\PP_2(D)$, and $F(\xi)=(\xi+\xi^2)/2$.
The map is bi-Lipschitz, whereas
$(x^2)\circ F=(\xi^2+2\xi^3+\xi^4)/4\notin\PP_2(D)$.

For each integer $\ell\ge1$, define
\begin{equation}\label{eq:intro-product}
W_h^\ell=\Span\left\{\prod_{\alpha=1}^{\ell}B_{i_\alpha,h}:
i_\alpha\in\mathcal I_h\right\},\qquad
V_h^\ell=\{w_h\circ F^{-1}:w_h\in W_h^\ell\}.
\end{equation}
The partition of unity gives $W_h^\ell\subseteq W_h^{\ell+1}$.
If $\PP_k(D)\subseteq W_h^1$, factorization of monomials gives
\begin{equation}\label{eq:intro-two-domain}
\PP_{\ell k}(D)\subseteq W_h^\ell,\qquad
F^*\PP_\ell(\Omega)\subseteq W_h^\ell,\qquad
\PP_\ell(\Omega)\subseteq V_h^\ell.
\end{equation}
The second inclusion follows independently from the geometry representation.
Neither inclusion requires independent product generators.

Our analytic construction starts from a locally $L^\infty$-stable operator
$Q_h:C(\overline D)\to W_h^1$ reproducing $\PP_k$.
Scaled local polynomial projections and a center-independent lifting yield
$\mathcal Q_h^\ell:L^2(D)\to W_h^\ell$, with exactness on $\PP_{\ell k}$
and local $L^2$ and $L^\infty$ stability. Under additional derivative
stability and global regularity of the initial space, the operator satisfies
\[
\|u-\mathcal Q_h^\ell u\|_{H^m(D)}
\le C_\ell h^{s-m}|u|_{H^s(D)},
\qquad 1\le s\le\ell k+1,
\]
for the integer derivative orders permitted by the hypotheses. The
constants are independent of $h$ for fixed $k$ and $\ell$.
Best-approximation estimates as $\ell$ increases are established separately
by polynomial containment; they do not assert uniform operator stability
in the product order.

For every fixed product order, the normalized generators yield an explicit
initial-space injection and a compatible redundant Galerkin system.
A local spectral coarse correction gives an energy spectral bound and
a computable stopping criterion linked to the approximation estimate.
The construction uses the GenEO principle~\cite{Spillane}; its
product-specific ingredient is the injection for arbitrary $\ell$.
The bound retains the support overlap explicitly and does not imply
uniform complexity as the product order increases.

Section~\ref{sec:products} states the assumptions and proves the algebraic
reproduction properties. Sections~\ref{sec:tools} and~\ref{sec:operator}
construct and analyze the quasi-interpolant. Section~\ref{sec:error}
establishes parametric and physical Sobolev estimates, and
Section~\ref{sec:pconv} treats increasing product order.
Section~\ref{sec:computation} develops the computational realization and
specifies the scope of its stability guarantees.

\section{Product spaces and simultaneous polynomial reproduction}
	\label{sec:products}
	\subsection{Mesh and initial-operator assumptions}\label{sec:setting}
	
For the common mesh scale $h\in(0,1]$, suppose $\mathcal T_h$ is a finite partition of $D$, up to sets of measure zero, into open elements $K$. 
\begin{assumption}\label{ass:mesh}
	There are constants $c_0,C_0>0$, independent of $h$ and $K$, such that
	\begin{equation}\label{eq:mesh}
		|K|\ge c_0h^d,\qquad \diam K\le C_0h
		\qquad(K\in\mathcal T_h).
	\end{equation}
	Here and below $|K|$ is Lebesgue measure. The diameter bound implies $|K|\le C_{\mathcal{T}} h^d$ by containment in a ball of radius $C_0h$, where
	$C_{\mathcal{T}}$ is independent of $h$ and $K$.
\end{assumption}

\begin{assumption}\label{ass:base}
	For a fixed integer $k\ge1$, a linear map
	$Q_h:C(\overline D)\to W_h^1$ satisfies
	\begin{align}
		Q_h(p|_{\overline D})&=p|_D,&& (p\in\PP_k),\label{eq:base-exact}\\
		\|Q_hv\|_{L^\infty(K)}&\le C_Q\|v\|_{L^\infty(\omega_K)},
		&&(v\in C(\overline D),\ K\in\mathcal T_h),\label{eq:base-stable}
	\end{align}
	where $\omega_K\subseteq D$ is open, $K\subseteq\omega_K$, and
	\begin{equation}\label{eq:patch}
		\diam\omega_K\le C_\omega h.
	\end{equation}
	The constants $C_Q,C_\omega>0$ are independent of $h$ and $K$.
\end{assumption}

The basic construction and all $L^2$ error estimates require only
Assumptions~\ref{ass:mesh}--\ref{ass:base} and \eqref{eq:partition_unity}.
For derivatives we use the following separate conditions.
\begin{assumption}
\label{ass:derivatives}
	For a fixed integer $M\ge1$, the initial space satisfies
	\begin{equation}\label{eq:base-regularity}
		W_h^1\subset W^{M,\infty}(D).
	\end{equation}
	For $1\le j\le M$ there is a constant $C_{Q,j}$, independent of $h,K,v$,
	such that
	\begin{equation}\label{eq:base-derivative}
		|Q_hv|_{H^j(K)}\le C_{Q,j}h^{d/2-j}\|v\|_{L^\infty(\omega_K)}
		,\quad(v\in C(\overline D),\ K\in\mathcal T_h).
	\end{equation}
\end{assumption}

\begin{proposition}
	\label{prop:base-inverse}
	Under Assumptions~\ref{ass:mesh}--\ref{ass:base}, if for $1\le j\le M$
	\begin{equation}\label{eq:base-inverse}
		|w_h|_{H^j(K)}\le C_{I,j}h^{-j}\|w_h\|_{L^2(K)},
		\quad(w_h\in W_h^1,\ K\in\mathcal T_h),
	\end{equation}
	where  $C_{I,j}>0$ is a constant
	independent of $h$, $K$, and $w_h$,
 then \eqref{eq:base-derivative} holds. 
\end{proposition}
\begin{proof}
The assumed inverse estimate \eqref{eq:base-inverse} applies with $w_h=Q_hv$, giving
	\begin{equation}\label{eq:base-inverse-proof-step1}
		|Q_hv|_{H^j(K)}
		\le C_{I,j}h^{-j}\|Q_hv\|_{L^2(K)}.
	\end{equation}
Since $Q_hv\in L^\infty(K)$ and $K$ has finite measure, we have
	\[
	\|Q_hv\|_{L^2(K)}^2
	=
	\int_K |Q_hv(\xi)|^2\,d\xi
	\le
	|K|\,\|Q_hv\|_{L^\infty(K)}^2.
	\]
	Then we obtain
	\begin{equation}\label{eq:base-inverse-proof-step2}
		\|Q_hv\|_{L^2(K)}
		\le |K|^{1/2}\|Q_hv\|_{L^\infty(K)}.
	\end{equation}
	Using the measure bound $|K|\le C_{\mathcal{T}}h^d$ and combining
	\eqref{eq:base-inverse-proof-step1}--
	\eqref{eq:base-inverse-proof-step2} with
	\eqref{eq:base-stable}, we obtain
	\begin{align*}
		|Q_hv|_{H^j(K)}
		&\le C_{I,j}h^{-j}|K|^{1/2}
		\|Q_hv\|_{L^\infty(K)}
		\\
		&\le C_{I,j}C_{\mathcal{T}}^{1/2}C_Q
		h^{d/2-j}\|v\|_{L^\infty(\omega_K)}.
	\end{align*}
	Thus \eqref{eq:base-derivative} holds with
	$C_{Q,j}=C_{I,j}C_{\mathcal{T}}^{1/2}C_Q$,
	which is independent of $h$, $K$, and $v$.
\end{proof}

\begin{remark}[A nonpolynomial way to verify the initial inverse estimate]
	\label{rem:reference}
	For example, suppose that each $K$ is a translate of $h\widehat K$, with
	$\widehat K$ chosen from finitely many bounded open reference sets, and
	that $w_h(a_K+h\,\cdot)|_{\widehat K}$ belongs to a fixed finite-dimensional
	space $E_{\widehat K}\subset H^M(\widehat K)$. On this space the $L^2$ norm
	is a norm and each $H^j$ seminorm is continuous, so
	$|e|_{H^j(\widehat K)}\le C\|e\|_{L^2(\widehat K)}$.
	Taking the largest constant over the finite collection and scaling yields
	\eqref{eq:base-inverse}. The spaces $E_{\widehat K}$ need not contain only
	polynomials. This is a sufficient verification mechanism, not an extra
	standing assumption.
\end{remark}

Unless explicitly stated otherwise, constants below may depend on the
fixed dimension, the indicated derivative and reproduction orders, and
the uniform constants in the invoked assumptions, but not on $h$, the
element, or the input function. We write $C_\ell$ when dependence on the
fixed product order is relevant. The local operator estimates are not
asserted to be uniform as $\ell\to\infty$. The best-approximation estimates
in Section~\ref{sec:pconv} have the separate constant dependence stated there.

\subsection{Generating functions and simultaneous polynomial reproduction}

For a multiplicity vector
${\kappa}
=(\kappa_1,\kappa_2,\ldots,\kappa_{N_h})\in\NN^{N_h}$, write
\(
|{\kappa}|
:=\sum_{i=1}^{N_h}\kappa_i.
\)
Fix an integer $\ell\ge1$. For $|\kappa|=\ell$, define the generating functions
\begin{equation}\label{eq:generators}
	G_{{\kappa},h}^{\ell}(\xi)
	:=
	\frac{\ell!}{\kappa_1!\kappa_2!\cdots\kappa_{N_h}!}
	\prod_{i=1}^{N_h}[B_{i,h}(\xi)]^{\kappa_i}.
\end{equation}
Here $\kappa_i$ records the multiplicity of $B_{i,h}$. Equivalently, the space in \eqref{eq:intro-product} satisfies
\[
W_h^\ell
=
\Span\left\{
G_{{\kappa},h}^{\ell}(\xi):
{\kappa}\in\NN^{N_h},\
|{\kappa}|=\ell
\right\}.
\]
For convenience set $W_h^0=\Span\{1\}$. All these spaces are finite-dimensional subspaces of $L^\infty(D)$.
The partition of unity property \eqref{eq:partition_unity} yields
\begin{equation}\label{eq:product-pou}
	1=\left(\sum_{i=1}^{N_h}B_{i,h}(\xi)\right)^\ell=\sum_{\substack{{\kappa}\in\NN^{N_h}\\
			|{\kappa}|=\ell}}\frac{\ell!}{\kappa_1!\kappa_2!\cdots\kappa_{N_h}!}
		\prod_{i=1}^{N_h}[B_{i,h}(\xi)]^{\kappa_i}=\sum_{\substack{{\kappa}\in\NN^{N_h}\\
			|{\kappa}|=\ell}}
	G_{{\kappa},h}^{\ell}(\xi),
\end{equation}
Each generator in \eqref{eq:generators} is nonnegative whenever all initial basis functions are nonnegative.

\begin{proposition}[Multiplicative structure and nestedness]
	\label{prop:algebra}
	For all integers $a,b\ge1$,
	\begin{equation}\label{eq:algebra}
		\Span\{uv:u\in W_h^a,\ v\in W_h^b\}=W_h^{a+b}.
	\end{equation}	
	In particular, the following nestedness  holds under
	the partition-of-unity condition \eqref{eq:partition_unity},
	\begin{equation}\label{eq:nestedness}
		W_h^\ell\subseteq W_h^{\ell+1}
		\qquad(\ell\ge0).
	\end{equation}
\end{proposition}
\begin{proof}
	
	For $u\in W_h^a$ and $v\in W_h^b$, the definition of the
	product spaces gives finite representations
	\[
	u=\sum_{r=1}^{R}c_r\prod_{i=1}^{a}u_{r,i},
	\qquad
	v=\sum_{s=1}^{S}d_s\prod_{j=1}^{b}v_{s,j},
	\]
	where $R,S$ are positive integers, $c_r,d_s\in\mathbb R$,
	and all $u_{r,i},v_{s,j}\in W_h^1$.
	Consequently,
	\[
	uv
	=
	\sum_{r=1}^{R}\sum_{s=1}^{S}c_rd_s
	\left(\prod_{i=1}^{a}u_{r,i}\right)
	\left(\prod_{j=1}^{b}v_{s,j}\right)
	\in W_h^{a+b}.
	\]
	Taking linear combinations proves
	\[
	\Span\{uv:u\in W_h^a,\ v\in W_h^b\}
	\subseteq W_h^{a+b}.
	\]
	Conversely, every generating product of $W_h^{a+b}$ can
	be written as
	\[
	\prod_{i=1}^{a+b}w_i
	=
	\left(\prod_{i=1}^{a}w_i\right)
	\left(\prod_{i=a+1}^{a+b}w_i\right),
	\qquad w_i\in W_h^1.
	\]
	The two factors belong to $W_h^a$ and $W_h^b$,
	respectively. Taking their linear span proves the reverse
	inclusion and hence \eqref{eq:algebra} for $a,b\ge1$.
	Finally, if the partition-of-unity property \eqref{eq:partition_unity} holds, that is $1\in W_h^1$, then for every $w_h\in W_h^\ell$,
	\[
	w_h=w_h\,1
	\in
	\Span\{uv:u\in W_h^\ell,\ v\in W_h^1\}
	=
	W_h^{\ell+1}.
	\]
	This proves \eqref{eq:nestedness} for $\ell\ge1$; the case
	$\ell=0$ follows from $1\in W_h^1$.
\end{proof}

\begin{theorem}[Simultaneous parametric and physical polynomial reproduction]
	\label{thm:reproduction}
	If \eqref{eq:base-exact} in Assumption \ref{ass:base} holds,
	then for every integer $\ell\ge1$, we have
\begin{equation}\label{eq:para-domain}
	\PP_{\ell k}(D)\subseteq W_h^\ell.
\end{equation}
	Independently, under the geometry representation \eqref{eq:ext_geo} and the partition of unity property \eqref{eq:partition_unity}, we have
	\begin{equation}\label{eq:phy-domain}
		F^*\PP_\ell(\Omega)\subseteq W_h^\ell.
	\end{equation}
\end{theorem}
\begin{proof}
	We first establish the polynomial product identity
	\begin{equation}\label{eq:polynomial-product}
		\PP_{p+q}(D)
		=
		\Span\{uv:u\in\PP_p(D),\ v\in\PP_q(D)\}
		\qquad(p,q\in\NN).
	\end{equation}
	The inclusion from right to left follows from the degree
	bound for a product of polynomials.
	For the reverse inclusion, consider a monomial $\xi^\gamma$
	with $|\gamma|\le p+q$.
	There exists a multi-index $\beta\le\gamma$ such that
	\[
	|\beta|=\min\{p,|\gamma|\}.
	\]
	Indeed, the prescribed total can be assigned successively
	to the components of $\beta$, without exceeding the
	corresponding components of $\gamma$, because it does not
	exceed $|\gamma|$.
	It follows that
	\[
	|\beta|\le p,
	\qquad
	|\gamma-\beta|
	=
	|\gamma|-\min\{p,|\gamma|\}
	\le q.
	\]
	Thus
	\[
	\xi^\gamma
	=
	\xi^\beta\xi^{\gamma-\beta},
	\qquad
	\xi^\beta\in\PP_p,
	\quad
	\xi^{\gamma-\beta}\in\PP_q.
	\]
	Restricting to $D$ and taking linear combinations proves
	\eqref{eq:polynomial-product}.
	
	We now prove the first inclusion in \eqref{eq:para-domain}
	by induction on $\ell$.
	The case $\ell=1$ is the assumed containment
	$\PP_k(D)\subseteq W_h^1$.
	If $\PP_{\ell k}(D)\subseteq W_h^\ell$, then
	\eqref{eq:polynomial-product} and
	Proposition~\ref{prop:algebra} give
	\[
	\begin{aligned}
		\PP_{(\ell+1)k}(D)
		&=
		\Span\{uv:u\in\PP_{\ell k}(D),\
		v\in\PP_k(D)\}\\
		&\subseteq
		\Span\{uv:u\in W_h^\ell,\ v\in W_h^1\}\\
		&=W_h^{\ell+1}.
	\end{aligned}
	\]
	This completes the induction and proves
	\begin{equation}\label{eq:param-repro}
		\PP_{\ell k}(D)\subseteq W_h^\ell.
	\end{equation}
	
	For physical reproduction, the geometry representation \eqref{eq:ext_geo} implies
	$F_\mu\in W_h^1$ for $1\le\mu\le d$.
	Moreover, $1\in W_h^1$ follows from
	the partition of unity property.
	For a physical monomial $x^\alpha$ with $|\alpha|\le\ell$,
	\[
	F^*(x^\alpha)
	=
	F_1^{\alpha_1}F_2^{\alpha_2}\cdots F_d^{\alpha_d}
	1^{\ell-|\alpha|}
	\in W_h^\ell.
	\]
	The expression on the right is a product of exactly $\ell$
	members of $W_h^1$, counting the constant factors.
	Since these monomials $x^\alpha$ with $|\alpha|\le\ell$ span $\PP_\ell(\Omega)$, linearity
	gives
	\[
	F^*\PP_\ell(\Omega)\subseteq W_h^\ell.
	\]
	This proves \eqref{eq:phy-domain}.
\end{proof}

These containments alone do not establish mesh-uniform
	approximation estimates.
	To obtain such estimates, the following sections construct
	a locally stable, polynomial-preserving operator with
	range in $W_h^\ell$.

\section{Local polynomial reconstruction and lifting}\label{sec:tools}

This section develops the auxiliary constructions needed for the
locally stable quasi-interpolant into \(W_h^\ell\) in Section~\ref{sec:operator}.
We introduce an auxiliary partition of unity and associated local
patches, define local polynomial reconstructions by \(L^2\)-projection,
and construct a bivariate polynomial lifting.
\subsection{Auxiliary partition of unity and local patches}

We first construct an auxiliary partition of unity using tensor products of one-dimensional hat functions on a uniform grid with mesh size comparable to \(h\). This construction is independent of the element geometry in \(\mathcal T_h\).

To this end, set $n_h=\lceil h^{-1}\rceil,~\delta_h=n_h^{-1}$, so that \(h/2\le\delta_h\le h\). For $j=0,1,\ldots,n_h$ and $t\in[0,1]$, define the one-dimensional hat functions by
\[
\vartheta_{j,h}(t)
=\max\{1-|t-j\delta_h|/\delta_h,0\}.
\]
Let $\mathcal A_h=\{0,1,\ldots,n_h\}^d$ be the index set of the auxiliary partition of unity.
For each \(a\in\mathcal A_h\) and $\xi=(\xi_1,\xi_2,\ldots,\xi_d)\in\overline D=[0,1]^d$, define the tensor-product function
\[
\theta_{a,h}(\xi)=\prod_{\nu=1}^d\vartheta_{a_\nu,h}(\xi_\nu).
\]
Then \(\delta_h a\in\overline D\) is the unique peak point of \(\theta_{a,h}(\xi)\), in particular, $\theta_{a,h}(\delta_h a)=1$.
To describe the support of \(\theta_{a,h}(\xi)\), set  $u_{a,\nu}^-=\max\{0,(a_\nu-1)\delta_h\}$,  $u_{a,\nu}^+=\min\{1,(a_\nu+1)\delta_h\}$ for $\nu=1,2,\ldots,d$, and define the associated open box \(U_{a,h}\) and its center \(z_{a,h}\) by
\[
U_{a,h}=\prod_{\nu=1}^d(u_{a,\nu}^-,u_{a,\nu}^+),
\qquad
(z_{a,h})_\nu=\frac{u_{a,\nu}^-+u_{a,\nu}^+}{2},
\quad \nu=1,2,\ldots,d.
\]
Then $\operatorname{supp}_{\overline D}\theta_{a,h}=\overline{U_{a,h}}$.
The functions \(\theta_{a,h}\) are continuous and nonnegative on
\(\overline D\), and satisfy
\begin{equation}\label{eq:hat-pou}
	\sum_{a\in\mathcal A_h}\theta_{a,h}(\xi)=1,
	\qquad\forall\xi\in\overline D.
\end{equation}
\begin{figure}[htbp]
	\centering
	\includegraphics[width=\textwidth]
	{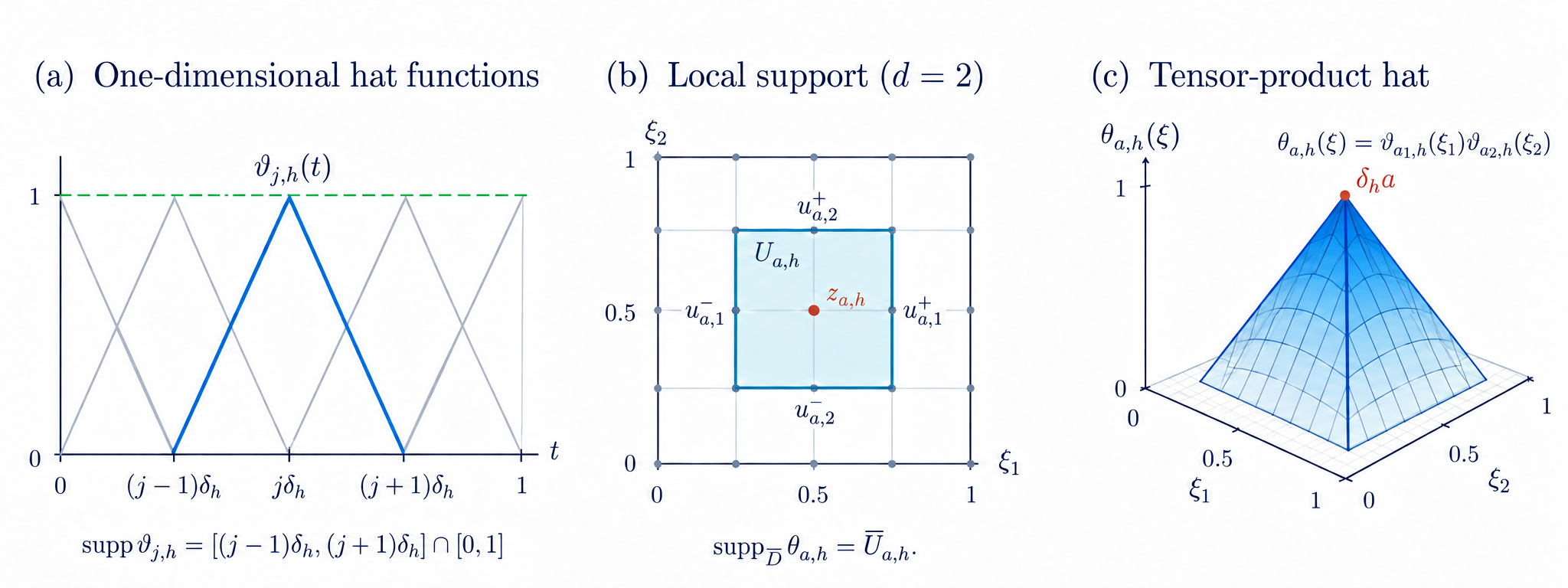}
	\caption{One-dimensional hat functions, a two-dimensional tensor-product hat function and the associated local support.}
	\label{fig:auxiliary-hat-functions}
\end{figure}

For each \(K\in\mathcal T_h\), let \(\omega_K\) be the open set appearing in Assumption~\ref{ass:base}.  Recall that \(K\subseteq\omega_K\) and
\(\operatorname{diam}\omega_K\le C_\omega h\).
We define the active index set $\mathcal A_K$ and the enlarged neighborhood $\omega_K^+$ by
\begin{equation}\label{eq:active-patches}
	\mathcal A_K:=\{a\in\mathcal{A}_h:\supp_{\overline D}\theta_{a,h}\cap\overline{\omega_K}
	\ne\varnothing\},
	\qquad
	\omega_K^+:=\bigcup_{a\in\mathcal A_K}U_{a,h}.
\end{equation}
The partition of unity implies $K\subseteq\omega_K\subseteq\omega_K^+$.
\begin{lemma}[Uniform auxiliary geometry]\label{lem:cover}
	For every \(a\in\mathcal A_h\),
	\begin{equation}\label{eq:aux-balls}
		B(z_{a,h},h/8)\subset U_{a,h}
		\subset B(z_{a,h},2\sqrt d\,h)\cap D.
	\end{equation}
	There is a constant \(C:=\max\left\{\bigl(4(C_\omega+\sqrt d)+1\bigr)^d,\,
	2(C_\omega+2\sqrt d)\right\}\), independent of \(h\) and \(K\), such that
	\[
	\#\mathcal A_K\le C,
	\qquad
	\operatorname{diam}(\omega_K^+)\le Ch,
	\]
	for every \(K\in\mathcal T_h\). Furthermore,
	\begin{equation}\label{eq:scaled-bounds}
		\sup_{\xi\in K}
		\left|\frac{\xi-z_{a,h}}h\right|\le C
		\quad(a\in\mathcal A_K),
		\qquad
		\sup_{\xi\in\overline{U_{a,h}}}
		\left|\frac{\xi-z_{a,h}}h\right|\le\sqrt d
		\quad(a\in\mathcal A_h).
	\end{equation}
\end{lemma}
\begin{proof}
 Every side length of \(U_{a,h}\) lies between \(\delta_h\) and
 \(2\delta_h\). Hence \(z_{a,h}\) is at distance at least
 \(\delta_h/2>h/4\ge h/8\) from every face of \(U_{a,h}\), while every
 point of \(\overline{U_{a,h}}\) is at distance at most
 \(\sqrt d\,\delta_h\le\sqrt d\,h<2\sqrt d\,h \) from \(z_{a,h}\).  This proves
	\eqref{eq:aux-balls} and the second estimate in
	\eqref{eq:scaled-bounds}.
	
	 For each \(a\in\mathcal A_K\), choose 
$t_a\in \overline{U_{a,h}}\cap\overline{\omega_K}$.
	By the assumed bound \(\operatorname{diam}(\overline{\omega_K})=\operatorname{diam}\omega_K\le C_\omega h\) given in Assumption~\ref{ass:base}, we have 
	\(
	|\xi-t_a|\le C_\omega h
	~\text{for every }\xi\in K.
	\)
	Since \(t_a\in\overline{U_{a,h}}\), it follows that
	\[
	|\xi-z_{a,h}|\le|\xi-t_{a}|+|t_{a}-z_{a,h}|\le(C_\omega+\sqrt d)h,\quad\forall \xi\in K,
	\]
	which proves the first estimate in \eqref{eq:scaled-bounds}, with \(C\ge C_\omega+\sqrt d\).
    Next, the definitions of \(U_{a,h}\) and \(\delta_h a\) imply
	\[
	|\delta_h a-t_a|\le\sqrt d\,\delta_h\le\sqrt d\,h.
	\]
	Fix a point \(\xi_K\in K\), we have
	\[
	|\delta_h a-\xi_K|
	\le |\delta_h a-t_a|+|t_a-\xi_K|
	\le(\sqrt d+C_\omega)h.
	\]
	Thus \(\forall~ a\in\mathcal A_K\), peak points \(\delta_h a\in \overline B\bigl(\xi_K,(C_\omega+\sqrt d)h\bigr)\).
	Since the lattice spacing satisfies \(\delta_h\ge h/2\), a lattice-point count in a coordinate
	cube containing the ball gives
	\[
	\#\mathcal A_K
	\le \bigl(4(C_\omega+\sqrt d)+1\bigr)^d.
	\]
	
	Finally, for every \(a\in\mathcal A_K\) and
	\(t_a'\in U_{a,h}\),  the preceding estimates give
	\[
	|t_a^{\prime}-\xi_K|
	\le |t_a^{\prime}-\delta_h a|+|\delta_h a-\xi_K|
	\le (2\sqrt d+C_\omega)h.
	\]
	Moreover, since \(\xi\in  K\subseteq\omega_K\), we have $|\xi-\xi_K|\le \operatorname{diam}(\omega_K)\le C_\omega h$. 
	Consequently, 
	\(
	\omega_K^+\subseteq\overline B\bigl(\xi_K,(C_\omega+2\sqrt d)h\bigr),
	\)
	and therefore
	\[
	\operatorname{diam}(\omega_K^+)
	\le 2(C_\omega+2\sqrt d)h,
	\]
	which completes the proof.
\end{proof}

To pass from local approximation estimates to global estimates, we introduce boxes that contain the enlarged patches and have uniformly bounded overlap.

For each \(K\in\mathcal T_h\), choose \(\xi_K\in K\). Fix
\(L>C_\omega+2\sqrt d+1\), and define
\begin{equation}\label{eq:boxes}
	V_K:=D\cap\prod_{\nu=1}^d
	((\xi_K)_\nu-Lh,(\xi_K)_\nu+Lh)=\prod_{\nu=1}^{d}
	\left(
	\max\{0,(\xi_K)_\nu-Lh\},
	\min\{1,(\xi_K)_\nu+Lh\}
	\right).
\end{equation}
\begin{lemma}[Approximation boxes and bounded overlap]
	\label{lem:overlap}
	For each \(K\in\mathcal T_h\), the set \(V_K\) is an open box
	containing \(\omega_K^+\), and each of its side lengths lies in
	\([h,2Lh]\).
	Moreover, for every \(\xi\in D\), set
	\(
	\mathcal T_h(\xi)
	:=\{K\in\mathcal T_h:\xi\in V_K\},
	\) and we have
	\begin{equation}\label{eq:overlap}
		\#\mathcal T_h(\xi)
		\le
		\frac{\omega_d}{c_0}
		\bigl(\sqrt d\,L+C_0\bigr)^d:=N_{\mathrm{ov}},
	\end{equation}
	where \(\omega_d:=|B(0,1)|\) is the volume of the unit ball in
	\(\mathbb R^d\). In particular, the overlap bound is independent
	of \(h\).
\end{lemma}
\begin{proof}
	By the estimates established in the proof of Lemma~\ref{lem:cover}, for $y\in\omega_K^+ $,
	\(
	|y-\xi_K|\le (C_\omega+2\sqrt d)h<Lh.
	\)
	Because \(\omega_K^+\subset D\), it follows from the definition
	of \(V_K\) that 
	\[
	\omega_K^+\subset V_K.
	\]
	For each coordinate direction, the length of the corresponding side of \(V_K\) is therefore
	at most \(2Lh\) and at least
	\[
	\min\{Lh,1\}\ge h,
	\]
	where we used \(L>1\) and \(0<h\le1\).
	
	To prove the overlap bound, fix \(\xi\in D\).
	If \(K\in\mathcal T_h(\xi)\), then the definition of \(V_K\)
	yields
	\(
	|\xi-\xi_K|<\sqrt d\,Lh.
	\)
	Since \(\xi_K\in K\) and
	\(\operatorname{diam}K\le C_0h\) by Assumption \ref{ass:mesh}, every \(y\in K\) satisfies
	\[
	|y-\xi|
	\le |y-\xi_K|+|\xi_K-\xi|
	<(\sqrt d\,L+C_0)h.
	\]
	Consequently, for each $K\in\mathcal T_h(\xi)$,
	\(
	K\subset
	B\bigl(\xi,(\sqrt d\,L+C_0)h\bigr).
	\)
	According to Assumption \ref{ass:mesh},
	the elements $K\in\mathcal{T}_h$ are pairwise disjoint up to sets of measure zero,
	and \(|K|\ge c_0h^d\). Hence
	\[
	\begin{aligned}
		\#\mathcal T_h(\xi)\,c_0h^d
		&\le
		\sum_{K\in\mathcal T_h(\xi)}|K|\\
		&\le
		\left|B\bigl(\xi,(\sqrt d\,L+C_0)h\bigr)\right|\\
		&=\omega_d(\sqrt d\,L+C_0)^d h^d.
	\end{aligned}
	\]
	Dividing by \(c_0h^d\) proves \eqref{eq:overlap}.
\end{proof}

\subsection{Scaled local polynomial reconstruction}

Fix an integer $\ell\ge 1$ for the remainder of the construction.
Set
\[
\Gamma_{\ell k}
:=
\left\{
\gamma\in\mathbb N_0^d:|\gamma|\le \ell k
\right\},
\qquad
y_{a,h}(\xi)
:=
\frac{\xi-z_{a,h}}{h},
\qquad
a\in\mathcal A_h.
\]
We represent each local $L^2$ projection in the scaled monomial
basis $\{y_{a,h}^{\gamma}\}_{\gamma\in\Gamma_{\ell k}}$.

For $a\in\mathcal A_h$ and $f\in L^2(D)$, define $P_{a,h}^{\ell k}f\in\mathbb P_{\ell k}(D)$ as the unique
polynomial satisfying
\begin{equation}\label{eq:projection}
	\int_{U_{a,h}}
	\bigl(P_{a,h}^{\ell k}f-f\bigr)p\,\mathrm d\xi=0
	\qquad\text{for all }p\in\mathbb P_{\ell k}(D).
\end{equation}
Equivalently, the restriction of $P_{a,h}^{\ell k}f$ to $U_{a,h}$ is the
$L^2(U_{a,h})$-orthogonal projection of $f|_{U_{a,h}}$ onto
$\mathbb P_{\ell k}(U_{a,h})$.
Writing
\begin{equation}\label{eq:projection-expansion}
	P_{a,h}^{\ell k}f
	=
	\sum_{\gamma\in\Gamma_{\ell k}}
	c_{a,h,\gamma}^{\ell}(f)y_{a,h}^{\gamma},
\end{equation}
and testing \eqref{eq:projection} with
$p=y_{a,h}^{\beta}$, $\beta\in\Gamma_{\ell k}$, yields the linear
system
\[
\mathbf M_{a,h}^{\ell k}\mathbf c_{a,h}^{\ell}(f)
=\mathbf b_{a,h}(f).
\]
Fix an ordering of $\Gamma_{\ell k}$ and  define the matrix $\mathbf M_{a,h}^{\ell k}$ and column vectors
\[
\mathbf c_{a,h}^{\ell}(f)
:=
\bigl(c_{a,h,\gamma}^{\ell}(f)\bigr)_{\gamma\in\Gamma_{\ell k}},
\qquad
\mathbf b_{a,h}(f)
:=
\bigl(b_{a,h,\beta}(f)\bigr)_{\beta\in\Gamma_{\ell k}},
\]
where
\begin{align*}
	[\mathbf M_{a,h}^{\ell k}]_{\beta\gamma}
	&:=
	\int_{U_{a,h}}
	y_{a,h}^{\beta+\gamma}\,\mathrm d\xi,
	\qquad \beta,\gamma\in\Gamma_{\ell k},\\
	b_{a,h,\beta}(f)
	&:=
	\int_{U_{a,h}}
	f(\xi)y_{a,h}(\xi)^\beta\,\mathrm d\xi,
	\qquad \beta\in\Gamma_{\ell k}.
\end{align*}
Since \eqref{eq:aux-balls} implies $|U_{a,h}|>0$, the scaled monomials $\{y_{a,h}^{\gamma}\}_{\gamma\in\Gamma_{\ell k}}$ are linearly independent in $L^2(U_{a,h})$. 
Thus $\mathbf M_{a,h}^{\ell k}$ is symmetric
positive definite and, in particular, invertible. Consequently,
\begin{equation}\label{eq:coefficients}
	c_{a,h,\gamma}^{\ell}(f)
	=
	\sum_{\beta\in\Gamma_{\ell k}}
	\bigl[
	(\mathbf M_{a,h}^{\ell k})^{-1}
	\bigr]_{\gamma\beta}
	b_{a,h,\beta}(f),
	\qquad
	\gamma\in\Gamma_{\ell k}.
\end{equation}
Then, for every \(p\in\mathbb P_q(D)\),
\begin{equation}\label{eq:local-reconstruction-reproduction}
	P_{a,h}\bigl(p|_{U_{a,h}}\bigr)
	=
	p|_{U_{a,h}}.
\end{equation}

\begin{lemma}\label{lem:coeff}
	For each $a\in\mathcal A_h$ and $\gamma\in\Gamma_{\ell k}$, the coefficient functional $c_{a,h,\gamma}^{\ell}(\cdot)$ is linear and depends only on $f|_{U_{a,h}}$. Moreover, there exists a constant $C_\ell>0$, such that
	\begin{align}
		|c_{a,h,\gamma}^{\ell}(f)|
		&\le
		C_\ell h^{-d/2}
		\|f\|_{L^2(U_{a,h})},
		&& f\in L^2(D),
		\label{eq:coeff-l2}\\
		|c_{a,h,\gamma}^{\ell}(f)|
		&\le
		C_\ell
		\|f\|_{L^\infty(U_{a,h})},
		&& f\in L^\infty(D).
		\label{eq:coeff-linf}
	\end{align}
After increasing $C_\ell$ by the fixed factor $\#\Gamma_{\ell k}$,
the same right-hand sides also bound
\[
\sum_{\gamma\in\Gamma_{\ell k}}|c_{a,h,\gamma}^{\ell}(f)|.
\]
The constant is independent of $a,h$, and $f$ for fixed $d,k,\ell$.
\end{lemma}
\begin{proof}
	Set $\widehat U_{a,h}:=\frac{U_{a,h}-z_{a,h}}{h}.$
	By \eqref{eq:aux-balls}, we have $B\left(0,\frac18\right)\subseteq\widehat U_{a,h}\subseteq B(0,2\sqrt d)$.
	Using the same ordering of $\Gamma_{\ell k}$ as for $\mathbf M_{a,h}^{\ell k}$, define the Gram matrix of the monomials $\{y^\gamma\}_{\gamma\in\Gamma_{\ell k}}$ by
	\[
	\mathbf G^{\ell k}
	:=
	\left(
	\int_{B(0,1/8)}
	y^{\beta+\gamma}\,\mathrm dy
	\right)_{\beta,\gamma\in\Gamma_{\ell k}}.
	\]
	The monomials $\{y^\gamma\}_{\gamma\in\Gamma_{\ell k}}$ are
	linearly independent in $L^2(B(0,1/8))$. Therefore, the Gram
	matrix $\mathbf G^{\ell k}$ is symmetric positive definite. For a real vector $v=(v_\gamma)_{\gamma\in\Gamma_{\ell k}}\in\mathbb R^{|\Gamma_{\ell k}|}$, let $p_v(y):=\sum_{\gamma\in\Gamma_{\ell k}}v_\gamma y^\gamma$. Setting $C_{\ell,1}:=\lambda_{\min}(\mathbf G^{\ell k})>0$,
	we obtain
	\[
	\int_{B(0,1/8)}|p_v(y)|^2\,\mathrm dy
	=
	v^{\mathsf T}\mathbf G^{\ell k}v
	\ge
	C_{\ell,1} \|v\|_2^2.
	\]
    Here, $\|\cdot\|_2$ denotes the Euclidean norm for vectors and the induced operator norm for matrices.  The matrix $\mathbf G^{\ell k}$ is independent of $a$ and $h$, and therefore so is $C_{\ell,1}$.
	Using the change of variables $\xi=z_{a,h}+hy$ and $B\left(0,\frac18\right)\subseteq\widehat U_{a,h}$, we obtain
	\[
	\begin{aligned}
		v^{\mathsf T}\mathbf M_{a,h}^{\ell k}v
		&=
		h^d
		\int_{\widehat U_{a,h}}|p_v(y)|^2\,\mathrm dy\\
		&\ge
		h^d
		\int_{B(0,\frac{1}{8})}|p_v(y)|^2\,\mathrm dy\ge
		C_{\ell,1} h^d\|v\|_2^2.
	\end{aligned}
	\]
	Therefore, $\lambda_{\min}(\mathbf M_{a,h}^{\ell k})\ge C_{\ell,1} h^d$ and 
	\begin{equation}\label{eq:M_inv}
		\bigl\|(\mathbf M_{a,h}^{\ell k})^{-1}\bigr\|_2 \le C_{\ell,1}^{-1}h^{-d},
	\end{equation}
	
	Next, the Cauchy--Schwarz inequality yields, for every
	$\beta\in\Gamma_{\ell k}$,
	\[
	|b_{a,h,\beta}(f)|
	\le
	\|f\|_{L^2(U_{a,h})}
	\|y_{a,h}^{\beta}\|_{L^2(U_{a,h})}.
	\]
	Since $\Gamma_{\ell k}$ is finite, the constant $C_{\ell,2}:=\max_{\beta\in\Gamma_{\ell k}}\int_{B(0,2\sqrt d)}|y^\beta|^2\,\mathrm dy$ is finite and independent of $a$ and $h$. By the same change of variables $\xi=z_{a,h}+hy$ and $\widehat U_{a,h}\subseteq B(0,2\sqrt d)$, for $\forall~\beta\in{\Gamma_{\ell k}}$, we have 
	\[
	\begin{aligned}
		\|y_{a,h}^{\beta}\|_{L^2(U_{a,h})}^2
		&=
		h^d
		\int_{\widehat U_{a,h}}|y^\beta|^2\,\mathrm dy\\
		&\le
		h^d
		\int_{B(0,2\sqrt d)}|y^\beta|^2\,\mathrm dy
		\le
		C_{\ell,2} h^d.
	\end{aligned}
	\]
	Thus we have $|b_{a,h,\beta}(f)|\le C_{\ell,2}^{1/2} h^{d/2}\|f\|_{L^2(U_{a,h})}$. Summing the above squares of the componentwise estimates over $\beta\in\Gamma_{\ell k}$, we obtain
	\[
	\begin{aligned}
		\|\mathbf b_{a,h}(f)\|_2^2=
		\sum_{\beta\in\Gamma_{\ell k}}
		|b_{a,h,\beta}(f)|^2\le
		C_{\ell,2} C_{\ell,3} h^d
		\|f\|_{L^2(U_{a,h})}^2,
	\end{aligned}
	\]
	where $C_{\ell,3}:=\#\Gamma_{\ell k}=\binom{d+\ell k}{d}$.
    Consequently,
	\begin{equation}\label{eq:b}
	\|\mathbf b_{a,h}(f)\|_2
	\le
	C_{\ell,2}C_{\ell,3} h^{d/2}
	\|f\|_{L^2(U_{a,h})}.
	\end{equation}
	Combining the estimates \eqref{eq:M_inv} and \eqref{eq:b} gives
	\[
	\begin{aligned}
		|c_{a,h,\gamma}^{\ell}(f)|
		&\le
		\|\mathbf c_{a,h}^{\ell}(f)\|_2\\
		&\le
		\bigl\|(\mathbf M_{a,h}^{\ell k})^{-1}\bigr\|_2
		\|\mathbf b_{a,h}(f)\|_2\le
		\frac{C_{\ell,2}C_{\ell,3}}{C_{\ell,1}} h^{-d/2}
		\|f\|_{L^2(U_{a,h})},
	\end{aligned}
	\]
	which proves \eqref{eq:coeff-l2}.
	Since $\widehat U_{a,h}\subseteq B(0,2\sqrt d)$, we have $|U_{a,h}|=h^d|\widehat U_{a,h}|\le C_d h^d$ where $C_{\ell,4}:=|B(0,2\sqrt d)|$.
    Thus, for $f\in L^\infty(D)$,
	\[
	\|f\|_{L^2(U_{a,h})}
	\le
	|U_{a,h}|^{1/2}\|f\|_{L^\infty(U_{a,h})}
	\le
	C_{\ell,4}^{1/2}h^{d/2}\|f\|_{L^\infty(U_{a,h})}.
	\]
	Substitution into \eqref{eq:coeff-l2} proves
	\eqref{eq:coeff-linf}.
	Therefore, both \eqref{eq:coeff-l2} and
	\eqref{eq:coeff-linf} hold with the constant
	\[
	C_\ell
	:=
	\frac{C_{\ell,2}C_{\ell,3}}{C_{\ell,1}}\max\{1,C_{\ell,4}^{1/2}\}.
	\]
\end{proof}

\subsection{Polynomial lifting}
\label{subsec:polynomial-lifting}

To construct a quasi-interpolant with values in the product space
$W_h^\ell$, we will use a representation of a polynomial of degree
at most $\ell k$ in which the degrees associated with two factors are
bounded by $k$ and $(\ell-1)k$, respectively. We obtain this
representation through a lifting to two independent vector variables.
The lifting preserves the original polynomial on the diagonal and
separates the degree bounds between the two groups of variables.

For $\alpha,\gamma\in\mathbb N_0^d$, the relation $\alpha\le\gamma$
means that $\alpha_\nu\le\gamma_\nu$ for every
$\nu=1,2,\ldots,d$. In this case, we write
\[
\binom{\gamma}{\alpha}
:=
\prod_{\nu=1}^d
\binom{\gamma_\nu}{\alpha_\nu}.
\]
For $n\in\mathbb N_0$ and $j\in\mathbb Z$, we adopt the convention $\binom nj=0$ for $j<0$ or $j>n$.

Let $\eta,\zeta\in\mathbb R^d$ be independent vector variables.
We define the polynomial space
\[
\mathcal B_{k,(\ell-1)k}
:=
\operatorname{span}
\left\{
\eta^\alpha\zeta^\beta:
\alpha,\beta\in\mathbb N_0^d,\ 
|\alpha|\le k,\ 
|\beta|\le {(\ell-1) k}
\right\}.
\]
Thus, the total degrees with respect to $\eta$ and $\zeta$ are
bounded by $k$ and ${(\ell-1) k}$, respectively.
For a polynomial $g(\tau)=\sum_{j=0}^{n}a_j\tau^j$ in the scalar variable $\tau\in\mathbb{R}$, we denote its coefficient of \(\tau^k\) by
\[
\operatorname{coef}_{\tau^k}(g):=a_k,
\]
where \(a_k:=0\) if \(k>n\). The coefficients \(a_j\) may themselves be polynomials in other variables.
Identifying each $p\in\mathbb P_{\ell k}(D)$ with its unique polynomial extension to $\mathbb R^d$, we define the lifting operator $\mathcal R_\ell$ by
\begin{equation}\label{eq:lift}
	(\mathcal R_\ell p)(\eta,\zeta)
	:=
	\binom {\ell k} k^{-1}
\operatorname{coef}_{\tau^k}
	\left[
	(1+\tau)^{\ell k}
	p\left(
	\frac{\tau\eta+\zeta}{1+\tau}
	\right)
	\right].
\end{equation}
The expression inside the coefficient extraction is understood after
cancellation of the denominator. Indeed, for every
$\gamma\in\Gamma_{\ell k}$,
\[
(1+\tau)^{\ell k}
\left(
\frac{\tau\eta+\zeta}{1+\tau}
\right)^\gamma
=
(1+\tau)^{{\ell k}-|\gamma|}
(\tau\eta+\zeta)^\gamma,
\]
which is a polynomial in $\tau$, $\eta$, and $\zeta$. Hence,
\eqref{eq:lift} is well defined.

For $\gamma\in\Gamma_{\ell k}$ and $\alpha\le\gamma$, define
\begin{equation}\label{eq:weights}
	w_{\gamma,\alpha}^{\ell}
	:=
	\frac{
		\displaystyle
		\binom{\gamma}{\alpha}
		\binom{{\ell k}-|\gamma|}{k-|\alpha|}
	}{
		\displaystyle
		\binom {\ell k}k
	}.
\end{equation}
\begin{lemma}
	\label{lem:lift}
   The lifting operator defined by
	\eqref{eq:lift} satisfies
	\begin{equation}\label{eq:lift-monomial}
		\mathcal R_\ell(\xi^\gamma)(\eta,\zeta)
		=
		\sum_{\alpha\le\gamma}
		w_{\gamma,\alpha}^{\ell}
		\eta^\alpha\zeta^{\gamma-\alpha}
	\end{equation}
   for each $\gamma\in\Gamma_{\ell k}$. Each nonzero weight $w_{\gamma,\alpha}^{\ell}$ satisfies
	$|\alpha|\le k$ and $|\gamma-\alpha|\le {(\ell-1) k}$. Thus for every $p\in\mathbb P_{\ell k}(D)$, we have $\mathcal R_\ell p\in\mathcal B_{k,(\ell-1)k}$. Furthermore,
	\begin{equation}\label{eq:lift-identities}
	(\mathcal R_\ell p)(\xi,\xi)=p(\xi),
\end{equation}
and $\sum_{\alpha\le\gamma} w_{\gamma,\alpha}^{\ell}=1,~\gamma\in\Gamma_{\ell k}$.
\end{lemma}

\begin{proof}
	We first apply \eqref{eq:lift} to the monomial
	$p(\xi)=\xi^\gamma$, where $\gamma\in\Gamma_{\ell k}$. By the
	multivariate binomial formula,
		\begin{align*}
		(1+\tau)^{\ell k}
		p\left(
		\frac{\tau\eta+\zeta}{1+\tau}
		\right)&=(1+\tau)^{\ell k}
		\left(
		\frac{\tau\eta+\zeta}{1+\tau}
		\right)^\gamma\\
		&=
		(1+\tau)^{{\ell k}-|\gamma|}
		(\tau\eta+\zeta)^\gamma\\
		&=
		\sum_{j=0}^{{\ell k}-|\gamma|}
		\binom{{\ell k}-|\gamma|}{j}\tau^j
		\sum_{\alpha\le\gamma}\binom{\gamma}{\alpha}\tau^{|\alpha|}\eta^\alpha\zeta^{\gamma-\alpha}\\
		&=	\sum_{\alpha\le\gamma}
		\sum_{j=0}^{{\ell k}-|\gamma|}
		\binom{\gamma}{\alpha}
		\binom{{\ell k}-|\gamma|}{j}
		\tau^{j+|\alpha|}
		\eta^\alpha\zeta^{\gamma-\alpha}.
	\end{align*}
	Consequently,
	\begin{align*}
		\operatorname{coef}_{\tau^k}
		\left[
		(1+\tau)^{\ell k}
		\left(
		\frac{\tau\eta+\zeta}{1+\tau}
		\right)^\gamma
		\right]
		&=
		\operatorname{coef}_{\tau^k}
		\left[
	\sum_{\alpha\le\gamma}
	\sum_{j=0}^{{\ell k}-|\gamma|}
	\binom{\gamma}{\alpha}
	\binom{{\ell k}-|\gamma|}{j}
	\tau^{j+|\alpha|}
	\eta^\alpha\zeta^{\gamma-\alpha}
		\right]
		\\
		&=
		\sum_{\alpha\le\gamma}
		\binom{\gamma}{\alpha}
		\binom{{\ell k}-|\gamma|}{k-|\alpha|}
		\eta^\alpha\zeta^{\gamma-\alpha}.
	\end{align*}
	It follows that
	\begin{equation}\label{eq:lift-monomial-proof}
		\mathcal R_\ell(\xi^\gamma)(\eta,\zeta)
		=
		\sum_{\alpha\le\gamma}
		w_{\gamma,\alpha}^{\ell}
		\eta^\alpha\zeta^{\gamma-\alpha}.
	\end{equation}
	If $w_{\gamma,\alpha}^{\ell}\ne0$, the convention for binomial coefficients implies $0\le k-|\alpha|\le \ell k-|\gamma|$.
	The first inequality gives $|\alpha|\le k$, whereas the second gives $	|\gamma-\alpha|=|\gamma|-|\alpha|\le {\ell k}-k={({\ell-1}) k}$.
	In particular, every nonzero term in \eqref{eq:lift-monomial} belongs to
	$\mathcal B_{k,({\ell-1}) k}$. Since \(\operatorname{coef}_{\tau^k}\) is linear, the definition
	\eqref{eq:lift} shows that \(\mathcal R_\ell\) is linear. Consequently, for each $p\in\mathbb P_{\ell k}(D)$,
	\[
	\mathcal R_\ell p\in\mathcal B_{k,({\ell-1}) k}.
	\]
	Set $\eta=\zeta=\xi$ in \eqref{eq:lift}, 
	then we have
	\begin{align*}
		(\mathcal R_\ell p)(\xi,\xi)
		&=
		\binom{\ell k}{k}^{-1}
		\operatorname{coef}_{\tau^k}
		\left((1+\tau)^{\ell k}p(\xi)\right)\\
		&=
		\binom{\ell k}{k}^{-1}
		\operatorname{coef}_{\tau^k}
		\left(
		\sum_{j=0}^{\ell k}
		\binom{\ell k}{j}p(\xi)\tau^j
		\right)\\
		&=
		\binom{\ell k}{k}^{-1}
		\binom{\ell k}{k}p(\xi)
		=p(\xi).
	\end{align*}
	 Applying the above identity to $p(\xi)=\xi^\gamma,~\gamma\in\Gamma_{\ell k}$ and using \eqref{eq:lift-monomial}, we find $	\left(\sum_{\alpha\le\gamma}w_{\gamma,\alpha}^{\ell} \right)\xi^\gamma=\xi^\gamma$ and this proves
	$\sum_{\alpha\le\gamma} w_{\gamma,\alpha}^{\ell}=1$.
\end{proof}
Moreover, if $p(\xi)=\sum_{\gamma\in\Gamma_{\ell k}}a_\gamma(\frac{\xi-c}{\sigma})^\gamma$
for $c\in\mathbb{R}^d$ and $\sigma>0$, then according to \eqref{eq:lift-monomial} and
\[
\frac{\frac{\tau\eta+\zeta}{1+\tau}-c}{\sigma}=\frac{\tau(\frac{\eta-c}{\sigma})+\frac{\zeta-c}{\sigma}}{1+\tau},
\]
we have
\begin{equation}\label{eq:center-independent}
	(\mathcal R_\ell p)(\eta,\zeta)
	=
	\sum_{\gamma\in\Gamma_{\ell k}}
	a_\gamma
	\sum_{\alpha\le\gamma}
	w_{\gamma,\alpha}^{\ell}
	\left(\frac{\eta-c}{\sigma}\right)^\alpha
	\left(\frac{\zeta-c}{\sigma}\right)^{\gamma-\alpha}.
\end{equation}

\begin{remark}
	When $\ell=2$, the lifting is symmetric, that is, $(\mathcal R_2p)(\eta,\zeta)=(\mathcal R_2p)(\zeta,\eta)$ for $p\in\mathbb P_{2k}(D)$.
	Indeed, \eqref{eq:weights} and the symmetry of binomial coefficients $	w_{\gamma,\alpha}^{2}=w_{\gamma,\gamma-\alpha}^{2}$.
	When $\ell=1$, the only nonzero term in \eqref{eq:lift-monomial} corresponds to $\alpha=\gamma$, and hence $(\mathcal R_1p)(\eta,\zeta)=p(\eta)$ for $p\in\mathbb P_k$.
	Thus the same definition covers both the initial-space case
	$\ell=1$ and all product orders $\ell\ge2$.
\end{remark}

\section{A stable product-space quasi-interpolant}\label{sec:operator}

In this section, we construct a quasi-interpolant \(\mathcal Q_h^\ell\) in the product space \(W_h^\ell\). The construction combines the initial operator \(Q_h\) with the local polynomial reconstruction and lifting introduced in Section \ref{sec:tools}. Under Assumptions \ref{ass:mesh}--\ref{ass:base}, the resulting operator is polynomially exact and locally stable. Additional derivative stability and Sobolev regularity follow under Assumption \ref{ass:derivatives}. These properties will be used in Section \ref{sec:error} to derive approximation estimates.

\subsection{Construction, polynomial exactness, and local stability}

For $f\in L^2(D)$, fix an integer $\ell\ge1$ and define a linear operator 
\begin{equation}\label{eq:operator}
	\mathcal Q_h^\ell f
	=\sum_{a\in\mathcal A_h}\sum_{\gamma\in\Gamma_{\ell k}}
	c_{a,h,\gamma}^{\ell}(f)
	\sum_{\substack{\alpha\le\gamma\\w_{\gamma,\alpha}^{\ell}\ne0}}
	w_{\gamma,\alpha}^{\ell}
	Q_h(\theta_{a,h}y_{a,h}^{\alpha})
	y_{a,h}^{\gamma-\alpha},
\end{equation}
Consider a summand for which \(w_{\gamma,\alpha}^{\ell}\ne0\), then Lemma~\ref{lem:lift} gives \(|\alpha|\le k$ and $|\gamma-\alpha|\le (\ell-1)k\).
By the definition of \(Q_h\), we have \(Q_h\!\left(\theta_{a,h}y_{a,h}^{\alpha}\right) \in W_h^1\).
Furthermore, $y_{a,h}^{\gamma-\alpha}\in\mathbb P_{(\ell-1)k}(D)$ and Theorem~\ref{thm:reproduction}  implies \(\mathbb P_{(\ell-1)k}(D)\subseteq W_h^{\ell-1}\).
 For $\ell=1$, the second factor is $1\in W_h^0$, so membership follows directly. For $\ell\ge2$, by \eqref{eq:algebra} given in Proposition~\ref{prop:algebra}, every nonzero summand in \eqref{eq:operator} belongs to \(W_h^\ell\). Therefore,  $\mathcal Q_h^\ell f\in W_h^\ell$. The coefficient functionals \(c_{a,h,\gamma}^{\ell}\) given by \eqref{eq:coefficients} are linear and all sums in \eqref{eq:operator} are finite, thus \(\mathcal Q_h^\ell\) is linear. Consequently,  \eqref{eq:operator} defines a linear operator
 \[
 \mathcal Q_h^\ell:L^2(D)\longrightarrow W_h^\ell.
 \]

\begin{theorem}[Polynomial exactness and local stability]\label{thm:operator}
	Assume the partition-of-unity property \eqref{eq:partition_unity} and
	Assumptions~\ref{ass:mesh}--\ref{ass:base}. For each $p\in\PP_{\ell k}(D)$,
	\begin{equation}\label{eq:exact}
		\mathcal Q_h^\ell(p)=p,
	\end{equation}
	Moreover, there are constants, independent of \(h\),
	\(K\), and \(f\), such that, for every \(f\in L^2(D)\) and
	\(K\in\mathcal T_h\),
	\begin{align}
		\|\mathcal Q_h^\ell f\|_{L^\infty(K)}
		&\le
		 C^{(\ell-1) k+1}C_Q C_\ell d^{k/2} h^{-d/2}
		\|f\|_{L^2(\omega_K^+)},
		\label{eq:product-qi-l2-linf}
		\\
		\|\mathcal Q_h^\ell f\|_{L^2(K)}
		&\le
		C_{\mathcal{T}}^{1/2} C^{(\ell-1) k+1}C_Q C_\ell d^{k/2}
		\|f\|_{L^2(\omega_K^+)}.
		\label{eq:product-qi-l2-stability}
	\end{align}
	If \(f\in L^\infty(D)\), then
	\begin{equation}
		\label{eq:product-qi-linf-stability}
		\|\mathcal Q_h^\ell f\|_{L^\infty(K)}
		\le
		C^{(\ell-1) k+1}C_Q C_\ell d^{k/2}\|f\|_{L^\infty(\omega_K^+)}.
	\end{equation}
	In particular, \((\mathcal Q_h^\ell f)|_K\) depends only on
	\(f|_{\omega_K^+}\).
\end{theorem}
\begin{proof}
Let $p\in\PP_{\ell k}$. Local projection reproduces $p$, and
\eqref{eq:center-independent} expresses its lifting in every local
coordinate system. Write the same lifted polynomial in unscaled variables as
\[
(\mathcal R_\ell p)(\eta,\zeta)
=\sum_{|\beta|\le(\ell-1)k}r_\beta(\eta)\zeta^\beta,
\qquad r_\beta\in\PP_k.
\]
The coefficients $r_\beta$ are independent of the patch index.
Applying $Q_h$ to the first variable after multiplication by
$\theta_{a,h}$, and evaluating on the diagonal, gives
\begin{align*}
\mathcal Q_h^\ell p(\xi)
&=\sum_{|\beta|\le(\ell-1)k}\xi^\beta
  \sum_{a\in\mathcal A_h}Q_h(\theta_{a,h}r_\beta)(\xi)\\
&=\sum_{|\beta|\le(\ell-1)k}\xi^\beta Q_h r_\beta(\xi)
=(\mathcal R_\ell p)(\xi,\xi)=p(\xi).
\end{align*}
Here we used linearity, the auxiliary partition of unity, and exactness
of $Q_h$ on $\PP_k$. All identities hold almost everywhere on $D$.

Next, we prove the local stability.
	Fix \(K\in\mathcal T_h\). We first show that only indices in
	\(\mathcal A_K\) contribute to \(\mathcal Q_h^\ell f\) on \(K\).
	If \(a\notin\mathcal A_K\), then the definition \eqref{eq:active-patches} of \(\mathcal A_K\) gives $\theta_{a,h}(\xi)=0,~\forall~\xi\in\omega_K$
	Hence, by the local stability of \(Q_h\),
	\[
	\left\|
	Q_h\!\left(\theta_{a,h}y_{a,h}^{\alpha}\right)
	\right\|_{L^\infty(K)}
	\le
	C_Q
	\left\|
	\theta_{a,h}y_{a,h}^{\alpha}
	\right\|_{L^\infty(\omega_K)}
	=0.
	\]
	Consequently, for each $a\notin\mathcal A_K$, $Q_h\!\left(\theta_{a,h}y_{a,h}^{\alpha}\right)$ vanishes almost everywhere on $K$. Then the restriction of \(\mathcal Q_h^\ell f\) to \(K\) is given by
	\begin{align}
		\mathcal Q_h^\ell f\big|_K
		={}&
		\sum_{a\in\mathcal A_K}
		\sum_{\gamma\in\Gamma_{\ell k}}
		c_{a,h,\gamma}^{\ell}(f)
		\sum_{\substack{\alpha\le\gamma\\
				w_{\gamma,\alpha}^{\ell}\ne0}}
		w_{\gamma,\alpha}^{\ell}
		Q_h\!\left(\theta_{a,h}y_{a,h}^{\alpha}\right)
		y_{a,h}^{\gamma-\alpha}\bigg|_K.
		\label{eq:product-qi-active-representation}
	\end{align}
	Consider \(a\in\mathcal A_K\). For \(w_{\gamma,\alpha}^{\ell}\ne0\), Lemma~\ref{lem:lift} gives \(|\alpha|\le k$ and $|\gamma-\alpha|\le (\ell-1)k\). Thus for every \(|\alpha|\le k\), the bound \eqref{eq:scaled-bounds} over $\overline{U_{a,h}}$, together with \(0\le\theta_{a,h}\le1\), gives
	\begin{equation}\label{eq:bounda-poly-local}
			\left\|
		\theta_{a,h}y_{a,h}^{\alpha}
		\right\|_{L^\infty(\omega_K)}
		\le d^{|\alpha|/2}
		\le d^{k/2}.
	\end{equation}
	The local stability of \(Q_h\) given by Assumption \ref{ass:base} therefore implies
	\begin{equation}
		\label{eq:product-qi-first-factor}
		\left\|
		Q_h\!\left(\theta_{a,h}y_{a,h}^{\alpha}\right)
		\right\|_{L^\infty(K)}\le C_Q\left\|\theta_{a,h}y_{a,h}^{\alpha}\right\|_{L^\infty(\omega_K)}
		\le C_Qd^{k/2}.
	\end{equation}
	Moreover, for each $a\in\mathcal A_K$ and $|\gamma-\alpha|\le (\ell-1) k$, the bound \eqref{eq:scaled-bounds} on \(K\) gives
	\begin{equation}
		\label{eq:product-qi-second-factor}
		\left\|
		y_{a,h}^{\gamma-\alpha}
		\right\|_{L^\infty(K)}
		\le C^{|\gamma-\alpha|}\le C^{(\ell-1) k}.
	\end{equation}
	It follows from
	\eqref{eq:product-qi-active-representation},
	\eqref{eq:product-qi-first-factor}, and
	\eqref{eq:product-qi-second-factor} that
	\begin{align*}
		\|\mathcal Q_h^\ell f\|_{L^\infty(K)}
		\le{}&
	C^{(\ell-1) k}C_Qd^{k/2}
		\sum_{a\in\mathcal A_K}
		\sum_{\gamma\in\Gamma_{\ell k}}
		\left|c_{a,h,\gamma}^{\ell}(f)\right|
		\sum_{\substack{\alpha\le\gamma\\
				w_{\gamma,\alpha}^{\ell}\ne0}}
		\left|w_{\gamma,\alpha}^{\ell}\right|.
	\end{align*}
	Using \eqref{eq:coeff-l2} and the uniform
	bound on \(\#\mathcal A_K\) given by Lemma \ref{lem:cover}, we obtain
	\[
	\|\mathcal Q_h^\ell f\|_{L^\infty(K)}
	\le
	 C^{(\ell-1) k+1}C_Q C_\ell d^{k/2} h^{-d/2}
	\max_{a\in\mathcal A_K}
	\|f\|_{L^2(U_{a,h})}
	\]
	whenever \(\mathcal A_K\ne\varnothing\). Since
	\(U_{a,h}\subseteq\omega_K^+\) for all \(a\in\mathcal A_K\),
	this proves \eqref{eq:product-qi-l2-linf}. If
	\(\mathcal A_K=\varnothing\), then
	\eqref{eq:product-qi-active-representation} shows that
	\(\mathcal Q_h^\ell f=0\) almost everywhere on \(K\), and the same
	estimate holds.
	The diameter bound in Assumption~\ref{ass:mesh} implies \(|K|\le C_{\mathcal{T}}h^d\).
	Thus, by \eqref{eq:product-qi-l2-linf},
	\[
	\begin{aligned}
		\|\mathcal Q_h^\ell f\|_{L^2(K)}
		&\le
		|K|^{1/2}\|\mathcal Q_h^\ell f\|_{L^\infty(K)}\le
		C_{\mathcal{T}}^{1/2} C^{(\ell-1) k+1}C_Q C_\ell d^{k/2} \|f\|_{L^2(\omega_K^+)}.
	\end{aligned}
	\]
	This proves \eqref{eq:product-qi-l2-stability}.
	If \(f\in L^\infty(D)\), the same argument, using the
	coefficient estimate \eqref{eq:coeff-linf} gives
	\[
	\|\mathcal Q_h^\ell f\|_{L^\infty(K)}
	\le
	 C^{(\ell-1) k+1}C_Q C_\ell d^{k/2}\|f\|_{L^\infty(\omega_K^+)}.
	\]
	This proves \eqref{eq:product-qi-linf-stability}.
	
	Finally, let \(f,g\in L^2(D)\) satisfy \(f=g\) almost everywhere
on \(\omega_K^+\). By linearity and \eqref{eq:product-qi-l2-stability},
	\[
	\begin{aligned}
		\|\mathcal Q_h^\ell f-\mathcal Q_h^\ell g\|_{L^2(K)}=\|\mathcal Q_h^\ell(f-g)\|_{L^2(K)}\le
		C_{\mathcal{T}}^{1/2} C^{(\ell-1) k+1}C_Q C_\ell d^{k/2}\|f-g\|_{L^2(\omega_K^+)}
		=0.
	\end{aligned}
	\]
	Hence \(\mathcal Q_h^\ell f=\mathcal Q_h^\ell g\) almost everywhere on \(K\).
	This proves that \((\mathcal Q_h^\ell f)|_K\) depends only on
	\(f|_{\omega_K^+}\).	
\end{proof}

\subsection{Sobolev regularity and derivative stability}

We next estimate derivatives of the constructed approximant. The
argument uses the derivative stability of the initial operator and
the scaling of the polynomial factors.

\begin{lemma}[Product regularity]
	\label{lem:product-space-regularity}
	Under \eqref{eq:base-regularity} in Assumption~\ref{ass:derivatives}, for each \(\ell\ge1\),
	\[
	W_h^\ell
	\subseteq W^{M,\infty}(D)
	\subseteq H^M(D).
	\]
\end{lemma}
\begin{proof}
	The weak Leibniz rule implies that \(W^{M,\infty}(D)\) is closed
	under multiplication. More precisely, if
	\(u,v\in W^{M,\infty}(D)\), then for every multi-index
	\(\mu\in\mathbb N_0^d\) with \(|\mu|\le M\),
	\[
	\partial^\mu(uv)
	=
	\sum_{\rho\le\mu}
	\binom{\mu}{\rho}
	(\partial^\rho u)
	(\partial^{\mu-\rho}v).
	\]
	 Each term on the right belongs to \(L^\infty(D)\). Since \(W_h^1\subseteq W^{M,\infty}(D)\), induction and the definition of \(W_h^\ell\) give the first inclusion. The second follows from the boundedness of \(D\).
\end{proof}

\begin{theorem}[Derivative stability of the quasi-interpolant]
	\label{thm:derivative}
	Assume the hypotheses of Theorem~\ref{thm:operator} and
	Assumption~\ref{ass:derivatives}. Fix an integer \(\ell\ge1\).
	Then, for every integer \(j\) with \(0\le j\le M\), there
	exists a constant \(C_{\ell,j}^{\mathrm{der}}>0\), independent
	of \(h\), \(K\), and \(f\), such that
	\begin{equation}\label{eq:derivative-stability}
		|\mathcal Q_h^\ell f|_{H^j(K)}
		\le
		C_{\ell,j}^{\mathrm{der}}h^{-j}
		\|f\|_{L^2(\omega_K^+)}
		\qquad
		(f\in L^2(D),\ K\in\mathcal T_h).
	\end{equation}
\end{theorem}
\begin{proof}
	For \(j=0\), the result is
	\eqref{eq:product-qi-l2-stability}. Let \(1\le j\le M\), and fix
	\(K\in\mathcal T_h\).
	For \(a\in \mathcal A_K\) and an index \(\alpha\) occurring in
	\eqref{eq:operator}, define $u_{a,\alpha}:=Q_h\!\left(\theta_{a,h}y_{a,h}^{\alpha}\right)$.
	As \eqref{eq:bounda-poly-local} in the proof of Theorem~\ref{thm:operator},
	\[
    \left\|
	\theta_{a,h}y_{a,h}^{\alpha}
	\right\|_{L^\infty(\omega_K)}
	\le d^{k/2}.
	\]
	For \(1\le|\rho|\le M\), the derivative stability condition in
	Assumption~\ref{ass:derivatives} gives
	\begin{align*}
		\|\partial^\rho u_{a,\alpha}\|_{L^2(K)}
		&\le
		|u_{a,\alpha}|_{H^{|\rho|}(K)}
		\le
		C_{Q,|\rho|}
		h^{d/2-|\rho|}
		\left\|
		\theta_{a,h}y_{a,h}^{\alpha}
		\right\|_{L^\infty(\omega_K)}
	\le
		C_{Q,|\rho|}d^{k/2}
		h^{d/2-|\rho|}.
		\label{eq:first-factor-positive-order}
	\end{align*}
	For \(\rho=0\), from the local \(L^\infty\) stability of \(Q_h\) given in \eqref{eq:base-stable} and \(|K|\le C_{\mathcal{T}}h^d\), we have 
	\[
		\|u_{a,\alpha}\|_{L^2(K)}\le |K|^{1/2}
		\|u_{a,\alpha}\|_{L^\infty(K)}\le C_{\mathcal{T}}^{1/2}C_Qd^{k/2} h^{d/2}.
	\]
	Define $C_{d,1}:=d^{k/2}\max\left\{C_{\mathcal{T}}^{1/2}C_Q,\,C_{Q,1},\ldots,C_{Q,M}\right\}.$
	Thus for $0\le|\rho|\le M$, we obtain
	\begin{equation}
		\label{eq:first-factor-derivative}
		\|\partial^\rho u_{a,\alpha}\|_{L^2(K)}
		\le
		C_{d,1}h^{d/2-|\rho|}.
	\end{equation}
	
	We now estimate the derivatives of the polynomial factor. Let
	\(\beta\in\mathbb N_0^d\) satisfy \(|\beta|\le (\ell-1)k\). Direct
	differentiation gives
	\[
	\partial^\lambda y_{a,h}^{\beta}
	=
	\begin{cases}
		\displaystyle
		h^{-|\lambda|}
		\frac{\beta!}{(\beta-\lambda)!}
		y_{a,h}^{\beta-\lambda},
		& \lambda\le\beta,\\[1.2ex]
		0,
		& \lambda\not\le\beta.
	\end{cases}
	\]
	Here $\beta!:=\prod_{\nu=1}^{d}\beta_\nu!$ and $	(\beta-\lambda)!:=\prod_{\nu=1}^{d}(\beta_\nu-\lambda_\nu)!$.
For $\lambda\le\beta$, $\frac{\beta!}{(\beta-\lambda)!}\le \max\{1,(\ell-1)k\}^{|\lambda|}$
and the scaled-coordinate bound \eqref{eq:scaled-bounds} implies
\[
\|y_{a,h}^{\beta-\lambda}\|_{L^\infty(K)}
\le C^{|\beta-\lambda|}\le C^{(\ell-1)k}.
\]
Hence, for every \(a\in\mathcal A_K\), \(|\beta|\le {(\ell-1)k}\), and \(|\lambda|\le M\),
\begin{equation}
	\label{eq:second-factor-derivative}
	\|\partial^\lambda y_{a,h}^{\beta}\|_{L^\infty(K)}
	\le
	C^{(\ell-1)k} \max\{1,(\ell-1)k\}^{|\lambda|}h^{-|\lambda|}.
\end{equation}
Let \(\mu\in\mathbb N_0^d\) satisfy \(|\mu|=j\). By the weak
product rule, \eqref{eq:first-factor-derivative}, and
\eqref{eq:second-factor-derivative},
\begin{align}
	\left\|
	\partial^\mu
	\left(u_{a,\alpha}y_{a,h}^{\beta}\right)
	\right\|_{L^2(K)}
	&\quad\le
	\sum_{\rho\le\mu}
	\binom{\mu}{\rho}
	\|\partial^\rho u_{a,\alpha}\|_{L^2(K)}
	\|\partial^{\mu-\rho}y_{a,h}^{\beta}\|_{L^\infty(K)}
	\nonumber\\
	&\quad\le
	C_{d,1}C^{(\ell-1)k}
	\sum_{\rho\le\mu}
	\binom{\mu}{\rho}
	\max\{1,(\ell-1)k\}^{|\mu-\rho|}
	h^{d/2-|\rho|-|\mu-\rho|}
	\nonumber\\
	&\quad\le
	C_{d,1}C^{(\ell-1)k}
	\left[1+\max\{1,(\ell-1)k\}\right]^j h^{d/2-j}.
	\label{eq:product-factor-derivative}
\end{align}
In the last step, we used $|\rho|+|\mu-\rho|=|\mu|=j$
and the multi-index binomial identity
\[
\sum_{\rho\le\mu}
\binom{\mu}{\rho}\max\{1,(\ell-1)k\}^{|\mu-\rho|}
=
\left[1+\max\{1,(\ell-1)k\}\right]^{|\mu|}.
\]
	
We now apply \eqref{eq:product-factor-derivative} to the finite sum \eqref{eq:operator}. For  $\gamma\in\Gamma_{\ell k}$ and every nonzero $w_{\gamma,\alpha}^{\ell}$ defined in \eqref{eq:weights}, Lemma~\ref{lem:lift} gives $|\alpha|\le k$, $|\gamma-\alpha|\le (\ell-1)k$ and $\sum_{\substack{\alpha\le\gamma}}w_{\gamma,\alpha}^{\ell}=1$. Therefore, for \(|\mu|=j\), the triangle inequality, \eqref{eq:coeff-l2}, and \eqref{eq:product-factor-derivative} imply
\begin{align*}
	\|\partial^\mu\mathcal Q_h^\ell f\|_{L^2(K)}
	\le{}&
	\sum_{a\in\mathcal A_K}
	\sum_{\gamma\in\Gamma_{\ell k}}
	|c_{a,h,\gamma}^{\ell}(f)|
	\sum_{\substack{\alpha\le\gamma\\
			w_{\gamma,\alpha}^{\ell}\ne0}}
	w_{\gamma,\alpha}^{\ell}
	\left\|
	\partial^\mu
	\left(
	u_{a,\alpha}y_{a,h}^{\gamma-\alpha}
	\right)
	\right\|_{L^2(K)}
	\\
	\le{}&
	C_{\ell}
	C_{d,1}C^{(\ell-1)k}\left[1+\max\{1,(\ell-1)k\}\right]^j
	h^{-j}
	\sum_{a\in\mathcal A_K}
	\sum_{\gamma\in\Gamma_{\ell k}}
	\|f\|_{L^2(U_{a,h})}.
\end{align*}
Since $U_{a,h}\subseteq\omega_K^+$ for each $a\in\mathcal A_K$ , $\#\mathcal A_K\le C$ from Lemma~\ref{lem:cover} and 
$C_{\ell,3}=\#\Gamma_{\ell k}=\binom{d+\ell k}{d}$ in the proof of Lemma~\ref{lem:coeff},
we obtain
\begin{equation}
	\label{eq:single-derivative-stability}
	\|\partial^\mu\mathcal Q_h^\ell f\|_{L^2(K)}
	\le
	C_{\ell}
	C_{d,1}
	C_{\ell,3}
	C^{(\ell-1)k+1}\left[1+\max\{1,(\ell-1)k\}\right]^j
	h^{-j}
	\|f\|_{L^2(\omega_K^+)}.
\end{equation}
For \(j\ge1\), the number of multi-indices \(\mu\in\mathbb N_0^d\) with \(|\mu|=j\) is $\binom{d+j-1}{j}$.
For \(j=0\), there is one such multi-index. Consequently,
\eqref{eq:single-derivative-stability} and the definition of
the \(H^j(K)\) seminorm prove
\eqref{eq:derivative-stability}, with
\[
C_{\ell,j}^{\mathrm{der}}
:=
\binom{d+j-1}{j}^{1/2}
C_{\ell}
C_{d,1}
C_{\ell,3}
C^{(\ell-1)k+1}\left[1+\max\{1,(\ell-1)k\}\right]^j,
\qquad 1\le j\le M.
\]
For $j=0$, use the constant in \eqref{eq:product-qi-l2-stability}.
The constant is independent of \(h\), \(K\), and \(f\).
\end{proof}

\section{Sobolev approximation estimates}
\label{sec:error}

In this section, we derive approximation estimates for the product-space quasi-interpolant constructed in Section~\ref{sec:operator}. We first establish local and global estimates on the parametric domain. We then transfer these estimates to the physical domain under suitable regularity assumptions on the geometry map and its inverse.

\subsection{Parametric-domain estimates}
\label{subsec:parametric-error}

We begin with a polynomial approximation result on the auxiliary boxes \(V_K\) introduced in Lemma~\ref{lem:overlap}. The result is a box version of the standard polynomial approximation estimates used in finite element analysis; see, for example, \cite{BrennerScott2008}. We include a proof to specify a constant that is uniform with respect to the box size and position.
\begin{lemma}[Polynomial approximation on comparable boxes]
	\label{lem:polynomial-approximation}
	
	Let \(d,s\in\mathbb N\) with \(d\ge1\) and \(s\ge1\), and let
	\(h>0\) and \(L\ge\tfrac12\). Let \(V\subset\mathbb R^d\) be an
	open axis-parallel box whose side lengths belong to
	\([h,2Lh]\).
	Then, for every \(u\in H^s(V)\), there exists a polynomial
	\(p\in\PP_{s-1}\) such that, for every integer \(0\le j\le s\),
	\begin{equation}
		\label{eq:polynomial-approximation}
		|u-p|_{H^j(V)}
		\le
		C_{\mathrm{box}}\,
		h^{s-j}|u|_{H^s(V)},
	\end{equation}
	where $C_{\mathrm{box}}:=(2L)^s d^{s/2}$.
\end{lemma}
\begin{proof}
	Let \(\widehat D:=(0,1)^d\). We first work on \(\widehat D\)
	and use the Sobolev seminorm convention
	\[
	|v|_{H^j(\widehat D)}^2
	:=
	\sum_{|\alpha|=j}
	\|\partial^\alpha v\|_{L^2(\widehat D)}^2.
	\]
	Given \(\widehat u\in H^s(\widehat D)\), we construct a
	polynomial \(\widehat p\in\PP_{s-1}\) such that, for each $|\alpha|\le s-1$,
	\begin{equation}
		\label{eq:moment-polynomial}
		\int_{\widehat D}
		\partial^\alpha(\widehat u-\widehat p)\,d\xi=0.
	\end{equation}
	The number of these equations equals \(\dim\PP_{s-1}\).
	To prove that the corresponding linear system is invertible,
	suppose that \(\widehat p\in\PP_{s-1}\) satisfies
	\[
	\int_{\widehat D}\partial^\alpha\widehat p\,d\xi=0
	\qquad
	(|\alpha|\le s-1).
	\]
	For \(|\alpha|=s-1\), the derivative
	\(\partial^\alpha\widehat p\) is constant and hence vanishes.
	Thus the homogeneous terms of degree \(s-1\) in
	\(\widehat p\) vanish. Repeating the argument for degrees
	\(s-2,\ldots,0\) gives \(\widehat p=0\). The system is therefore
	injective and hence invertible.
	
	We next establish the zero-mean estimate, that is, for $v\in H^1(\widehat D)$ satisfying $\int_{\widehat D}v\,d\xi=0$,
	\begin{equation}
		\label{eq:unit-cube-poincare}
		\|v\|_{L^2(\widehat D)}
		\le
		|v|_{H^1(\widehat D)}.
	\end{equation}
	For \(i=1,2,\ldots,d\), let  $(A_i v)(\xi_1,\ldots,\xi_d):=
	\int_0^1
	v(\xi_1,\ldots,\xi_{i-1},t,\xi_{i+1},\ldots,\xi_d)
	\,dt$ denote averaging with respect
	to the \(i\)-th coordinate, and put
	\[
	P_0:=I,
	\qquad
	P_i:=A_1\cdots A_i.
	\]
	The operators \(A_i\) are commuting orthogonal projections in \(L^2(\widehat D)\). Consequently,
	\[
	v-P_dv
	=
	\sum_{i=1}^d(P_{i-1}v-P_iv)
	\]
	is an orthogonal decomposition in \(L^2(\widehat D)\).
	Since \(v\) satisfies $\int_{\widehat D}v\,d\xi=0$, we have \(P_dv=0\). The
	\(L^2\)-contractivity of \(P_{i-1}\) and the one-dimensional
	Poincar\'e inequality give
	\[
	\begin{aligned}
		\|v\|_{L^2(\widehat D)}^2=
		\sum_{i=1}^d
		\|P_{i-1}(v-A_iv)\|_{L^2(\widehat D)}^2\le
		\sum_{i=1}^d
		\|v-A_iv\|_{L^2(\widehat D)}^2
	\le
		\sum_{i=1}^d
		\|\partial_i v\|_{L^2(\widehat D)}^2
		=
		|v|_{H^1(\widehat D)}^2.
	\end{aligned}
	\]
	This proves \eqref{eq:unit-cube-poincare}.
	
	Put \(\widehat e:=\widehat u-\widehat p\). For every \(\alpha\) with \(|\alpha|\leq s-1\),
	\eqref{eq:moment-polynomial} and
	\eqref{eq:unit-cube-poincare} imply
	\[
	\|\partial^\alpha(\widehat u-\widehat p)\|_{L^2(\widehat D)}^2
	\le
	\sum_{i=1}^d
	\|\partial_i\left[\partial^\alpha(\widehat u-\widehat p)\right]\|_{L^2(\widehat D)}^2	=
	\left|
	\partial^\alpha(\widehat u-\widehat p)
	\right|_{H^1(\widehat D)}^2.
	\]
	Summing over \(|\alpha|=j\), and observing that each derivative
	of total order \(j+1\) occurs at most \(d\) times, we obtain
	\[
	|\widehat u-\widehat p|_{H^j(\widehat D)}^2=\sum_{|\alpha|=j}\|\partial^\alpha(\widehat u-\widehat p)\|_{L^2(\widehat D)}^2
	\le 
	d\,|\widehat u-\widehat p|_{H^{j+1}(\widehat D)}^2,
	\qquad
	0\le j<s.
	\]
	Since \(\widehat p\in\PP_{s-1}\), its derivatives of order
	\(s\) vanish. Iteration therefore yields
	\begin{equation}
		\label{eq:reference-poly-bound}
		|\widehat u-\widehat p|_{H^j(\widehat D)}
		\le
		d^{(s-j)/2}
		|\widehat u|_{H^s(\widehat D)},
		\qquad
		0\le j\le s.
	\end{equation}
	
	Write
	\[
	V=b+A\widehat D,
	\qquad
	A=\operatorname{diag}(a_1,\ldots,a_d),
	\qquad
	h\le a_i\le2Lh.
	\]
	Define
	\[
	\widehat u(\xi):=u(b+A\xi),
	\qquad
	p(x):=\widehat p\bigl(A^{-1}(x-b)\bigr).
	\]
	In particular, \(p\in\PP_{s-1}\). By the chain rule for the
	affine map and a change of variables, for \(0\le j\le s\),
	\[
	|u-p|_{H^j(V)}
	\le
	(\det A)^{1/2}h^{-j}
	|\widehat u-\widehat p|_{H^j(\widehat D)}
	\]
	and
	\[
	|\widehat u|_{H^s(\widehat D)}
	\le
	(\det A)^{-1/2}(2Lh)^s
	|u|_{H^s(V)}.
	\]
	Combining these inequalities with
	\eqref{eq:reference-poly-bound} gives
	\[
	|u-p|_{H^j(V)}
	\le
	(2L)^s d^{(s-j)/2}
	h^{s-j}|u|_{H^s(V)}.
	\]
	Since \(d\ge1\) and \(0\le j\le s\), we have $d^{(s-j)/2}\le d^{s/2}$.
	Thus \eqref{eq:polynomial-approximation} holds with the constant
	\[
	C_{\mathrm{box}}=(2L)^s d^{s/2}.
	\]
\end{proof}

According to Lemma~\ref{lem:overlap}, we have the following uniform bound for the auxiliary boxes $V_K$ defined in \eqref{eq:boxes},
\begin{equation*}
	\#\mathcal T_h(\xi)
	\le
	\frac{\omega_d}{c_0}
	\bigl(\sqrt d\,L+C_0\bigr)^d:=N_{\mathrm{ov}}.
\end{equation*}

For a fixed \(\ell\), with $C_{\ell,0}^{\mathrm{der}}:=C_{\mathcal{T}}^{1/2} C^{(\ell-1) k+1}C_Q C_\ell d^{k/2}$, Theorem~\ref{thm:operator} gives the local \(L^2\)-stability estimate
\begin{equation}
	\label{eq:parametric-stability-constant-zero}
	\|\mathcal Q_h^\ell f\|_{L^2(K)}
	\le
	C_{\ell,0}^{\mathrm{der}}\|f\|_{L^2(\omega_K^+)}.
\end{equation}
When Assumption~\ref{ass:derivatives} holds, for $1\le j\le M$, Theorem~\ref{thm:derivative} gives
\begin{equation}
	\label{eq:parametric-stability-constant-j}
	|\mathcal Q_h^\ell f|_{H^j(K)}
	\le
	C_{\ell,j}^{\mathrm{der}}h^{-j}
	\|f\|_{L^2(\omega_K^+)}.
\end{equation}

For the fixed \(s,m,\ell\), define
\begin{equation}
	\label{eq:parametric-error-constants}
	C_{\mathrm{st},\ell,m}^{\max}
	:=
	\max_{0\le j\le m}	C_{\ell,j}^{\mathrm{der}},
	\qquad
	C_{\mathrm{loc}}
	:=\bigl(1+C_{\mathrm{st},\ell,m}^{\max}\bigr)C_{\mathrm{box}},
\end{equation}
and
\begin{equation}
	\label{eq:global-error-constant}
	C_{\mathrm{par}}
	:=
	\sqrt{(m+1)N_{\mathrm{ov}}}\,C_{\mathrm{loc}}.
\end{equation}
When \(m=0\), the maximum in
\eqref{eq:parametric-error-constants} consists only of
\(C_{\ell,0}^{\mathrm{der}}\), and no derivative stability constant is
needed.

\begin{theorem}[Parametric-domain approximation]
\label{thm:error}
Assume the partition-of-unity property and Assumptions~\ref{ass:mesh}
and~\ref{ass:base}. Fix integers $\ell\ge1$, $1\le s\le\ell k+1$,
and $0\le m\le s$. If $m\ge1$, assume additionally
Assumption~\ref{ass:derivatives} with $M\ge m$.
Then, for $u\in H^s(D)$, $K\in\mathcal T_h$, and $0\le j\le m$,
\begin{equation}
		\label{eq:local-error}
		|u-\mathcal Q_h^\ell u|_{H^j(K)}
		\le
		C_{\mathrm{loc}}h^{s-j}|u|_{H^s(V_K)},
	\end{equation}
	and
	\begin{equation}
		\label{eq:global-error}
		\|u-\mathcal Q_h^\ell u\|_{H^m(D)}
		\le
		C_{\mathrm{par}}h^{s-m}|u|_{H^s(D)}.
	\end{equation}
	Consequently,
	\begin{equation}
		\label{eq:best-error}
		\inf_{w_h\in W_h^\ell}
		\|u-w_h\|_{H^m(D)}
		\le
		C_{\mathrm{par}}h^{s-m}|u|_{H^s(D)}.
	\end{equation}
	The constants \(C_{\mathrm{loc}}\) and \(C_{\mathrm{par}}\) are
	defined by \eqref{eq:parametric-error-constants} and
	\eqref{eq:global-error-constant}, respectively.
\end{theorem}

\begin{proof}
	Fix \(K\in\mathcal T_h\). By Lemma~\ref{lem:overlap}, the box \(V_K\)
	has side lengths in \([h,2Lh]\) and contains \(\omega_K^+\).
	Lemma~\ref{lem:polynomial-approximation} therefore provides a
	polynomial \(p_K\in\PP_{s-1}\) such that, for $0\le j\le s$,
	\begin{equation}
		\label{eq:local-polynomial-error}
		|u-p_K|_{H^j(V_K)}
		\le
		C_{\mathrm{box}}h^{s-j}
		|u|_{H^s(V_K)}.
	\end{equation}
	Since \(s-1\le\ell k\), the polynomial exactness of \(\mathcal Q_h^\ell\) in Theorem~\ref{thm:operator} gives $	\mathcal Q_h^\ell p_K=p_K$.
	Consequently,
	\[
	u-\mathcal Q_h^\ell u
	=
	(u-p_K)-\mathcal Q_h^\ell(u-p_K).
	\]
	For \(0\le j\le m\), the local stability estimates
	\eqref{eq:parametric-stability-constant-zero} and
	\eqref{eq:parametric-stability-constant-j} yield
	\[
	|\mathcal Q_h^\ell(u-p_K)|_{H^j(K)}
	\le
	C_{\mathrm{st},\ell,m}^{\max}h^{-j}
	\|u-p_K\|_{L^2(\omega_K^+)}.
	\]
	Since \(\omega_K^+\subseteq V_K\),
	\eqref{eq:local-polynomial-error} with \(j=0\) gives
	\[
	|\mathcal Q_h^\ell(u-p_K)|_{H^j(K)}
	\le
	C_{\mathrm{st},\ell,m}^{\max}
	C_{\mathrm{box}}
	h^{s-j}|u|_{H^s(V_K)}.
	\]
	Moreover, \(K\subseteq V_K\), so
	\[
	|u-p_K|_{H^j(K)}
	\le
	C_{\mathrm{box}}
	h^{s-j}|u|_{H^s(V_K)}.
	\]
	The triangle inequality proves \eqref{eq:local-error}.
	
	If \(m\ge1\), Lemma~\ref{lem:product-space-regularity} and 
	Assumption~\ref{ass:derivatives} imply
	\(\mathcal Q_h^\ell u\in H^m(D)\). For \(m=0\), the same conclusion
	holds with \(H^0(D)=L^2(D)\) by the mapping property of
	\(\mathcal Q_h^\ell\). Since the elements are pairwise disjoint and
	cover \(D\) up to a set of measure zero,
	\[
	\|u-\mathcal Q_h^\ell u\|_{H^m(D)}^2
	=
	\sum_{j=0}^{m}
	\sum_{K\in\mathcal T_h}
	|u-\mathcal Q_h^\ell u|_{H^j(K)}^2.
	\]
	Applying \eqref{eq:local-error} and the overlap bound \eqref{eq:overlap} yields
	\[
	\|u-\mathcal Q_h^\ell u\|_{H^m(D)}^2
	\le
	N_{\mathrm{ov}}C_{\mathrm{loc}}^2
	\sum_{j=0}^{m}h^{2(s-j)}
	|u|_{H^s(D)}^2.
	\]
	Because \(0<h\le1\) and \(j\le m\), we have $h^{2(s-j)}\le h^{2(s-m)}$.
	Thus
	\[
	\|u-\mathcal Q_h^\ell u\|_{H^m(D)}^2
	\le
	(m+1)N_{\mathrm{ov}}C_{\mathrm{loc}}^2
	h^{2(s-m)}|u|_{H^s(D)}^2,
	\]
	which proves \eqref{eq:global-error} with $C_{\mathrm{par}}=\sqrt{(m+1)N_{\mathrm{ov}}}\,C_{\mathrm{loc}}$.
	
	Finally, \(\mathcal Q_h^\ell u\in W_h^\ell\), so it is an admissible
	competitor in the infimum in \eqref{eq:best-error}.
\end{proof}

\begin{remark}[Interpretation of the rate]\label{rem:rate}
	At $\ell=1$, Theorem~\ref{thm:error} guarantees $1\le s\le k+1$. Product
	enrichment extends this range to $\ell k+1$. It does not assert that
	the initial space has exactly approximation order $k+1$, or that the
	displayed estimate is a lower bound on the error. In particular, an
	initial space with additional reproduction may have a better baseline.
	The available global derivative orders remain constrained by $M$:
	product enrichment does not automatically increase interelement
	regularity.
\end{remark}

\begin{remark}[Scope of the assumptions]\label{rem:scope}
	No local degree bound or polynomial structure of $W_h^1$ was used in
	Theorem~\ref{thm:error}. The local projections and the polynomial
	lifting are auxiliary analytic constructions. They do not identify
	$W_h^1$ with a polynomial space. Derivative stability
	\eqref{eq:base-derivative} is essential input to the derivative argument;
	its verification cannot be replaced by finite dimensionality at each
	fixed mesh level without uniform control.
	
	The mesh may have tensor-product, simplicial, or T-mesh connectivity
	provided the stated assumptions are verified. The common scale and
	volume bounds do not cover arbitrary grading with unbounded size
	ratios. General parameter domains require suitable auxiliary coverings
	and uniformly controlled polynomial-approximation patches. The
	results concern standard Sobolev spaces, not general broken or bent
	Sobolev spaces.
\end{remark}

\subsection{Physical-domain estimates}
\label{sec:physical}
We now transfer the parametric-domain estimates to the physical
domain. The geometry map \(F:D\to\Omega\) is the fixed bi-Lipschitz bijection in \eqref{eq:ext_geo}.
We first state the composition estimate used to transport the error. Let \(U,~V\subset\mathbb R^d\) be open bounded sets.
For a bi-Lipschitz bijection \(G:U\to V\), define
\begin{equation}
	\label{eq:jacobian-composition-constant}
	J_G
	:=
	\left(
	\operatorname*{ess\,sup}_{y\in V}
	|\det DG^{-1}(y)|
	\right)^{1/2}.
\end{equation}
For each $v\in L^2(V)$, the change-of-variables formula for bi-Lipschitz mappings \cite{EvansGariepy2015} gives
\begin{equation}
	\label{eq:l2-composition}
	\|v\circ G\|_{L^2(U)}
	\le
	J_G\|v\|_{L^2(V)}.
\end{equation}
Moreover, $J_G\le\bigl(\operatorname{Lip}(G^{-1})\bigr)^{d/2}<\infty$, where $\operatorname{Lip}(G^{-1})$ is the Lipschitz constant of \(G^{-1}\).

\begin{lemma}[Integer-order Sobolev composition]
	\label{lem:composition}
	
	Let $U,V$ be bounded open domains, let $G:U\to V$ be bi-Lipschitz, and let \(r\ge1\) be an integer. Suppose that
	\(G\in W^{r,\infty}(U;\mathbb R^d)\).
	Then for $v\in H^r(V)$, there exists a finite constant \(C_{\mathrm{comp}}(G,r)>0\) such that
	\begin{equation}
		\label{eq:composition}
		\|v\circ G\|_{H^r(U)}
		\le
		C_{\mathrm{comp}}(G,r)\|v\|_{H^r(V)}.
	\end{equation}
	The constant \(C_{\mathrm{comp}}(G,r)\) can be chosen to depend only on
	\(d\), \(r\), \(J_G\), and the \(L^\infty\) norms of the weak
	derivatives of \(G\) of orders \(1,2,\ldots,r\).
\end{lemma}
\begin{proof}
	Let \(v\in H^1(V)\). {By \cite[Theorem~2.46]{HajlaszEllipticNotes}, the weak chain rule for a bi-Lipschitz
	mapping gives}
	\begin{equation}
		\label{eq:composition-proof-chain}
		\partial_i(v\circ G)
		=
		\sum_{a=1}^d
		\bigl((\partial_a v)\circ G\bigr)\,\partial_iG_a,
		\qquad i=1,2,\ldots,d,
	\end{equation}
	where the derivatives are understood in the weak sense. In particular, \(v\circ G\in H^1(U)\).
	Set $L_1:=
	\operatorname*{ess\,sup}_{\xi\in U}
	\|DG(\xi)\|_{2}$.
	Using the vector form
	of \eqref{eq:composition-proof-chain} and \eqref{eq:l2-composition}, we obtain
	\[
		|v\circ G|_{H^1(U)}\le
		L_1\|(\nabla v)\circ G\|_{L^2(U)}\le
		J_G L_1 |v|_{H^1(V)}.
	\]
	Together with \eqref{eq:l2-composition}, this gives
	\begin{equation}
		\label{eq:composition-proof-h1}
		\|v\circ G\|_{H^1(U)}
		\le
		J_G\max\{1,L_1\}\|v\|_{H^1(V)}.
	\end{equation}

	Next, for every integer \(j\ge0\), the weak Leibniz rule yields
	a constant \(C_{\mathrm{mult}}(d,j)>0\), depending only on
	\(d\) and \(j\), such that
	\begin{equation}
		\label{eq:composition-proof-multiplier}
		\|aw\|_{H^j(U)}
		\le
		C_{\mathrm{mult}}(d,j)
		\|a\|_{W^{j,\infty}(U)}
		\|w\|_{H^j(U)}
	\end{equation}
	for all \(a\in W^{j,\infty}(U)\) and \(w\in H^j(U)\).
	
	We now prove the composition estimate at higher orders by
	induction. Suppose that, for all $w\in H^{n-1}(V)$ and  some integer \(n\) with
	\(2\le n\le r\), there is a finite constant \(C_{\mathrm{comp}}(G,n-1)\) such
	that
	\begin{equation}
		\label{eq:composition-proof-induction}
		\|w\circ G\|_{H^{n-1}(U)}
		\le
		C_{\mathrm{comp}}(G,n-1)\|w\|_{H^{n-1}(V)}.
	\end{equation}
	Define $A_{n-1}:=\max_{1\le i,a\le d}\|\partial_iG_a\|_{W^{n-1,\infty}(U)}$.
	This quantity is finite because
	\(G\in W^{r,\infty}(U;\mathbb R^d)\) and \(n\le r\).
	Let \(v\in H^n(V)\). For each \(a=1,2,\ldots,d\),
	\(\partial_a v\in H^{n-1}(V)\). Hence
	\eqref{eq:composition-proof-induction} gives
	\[
	\|(\partial_a v)\circ G\|_{H^{n-1}(U)} \le C_{\mathrm{comp}}(G,n-1)\|\partial_a v\|_{H^{n-1}(V)}.
	\]
	Applying \eqref{eq:composition-proof-multiplier} to each term
	in \eqref{eq:composition-proof-chain}, we find that, for
	\(i=1,2,\ldots,d\),
	\[
	\begin{aligned}
		\|\partial_i(v\circ G)\|_{H^{n-1}(U)}
		&\le
		C_{\mathrm{mult}}(d,n-1)
		\sum_{a=1}^d
		\|\partial_iG_a\|_{W^{n-1,\infty}(U)}
		\|(\partial_a v)\circ G\|_{H^{n-1}(U)}
		\\
		&\le
		C_{\mathrm{mult}}(d,n-1)
		A_{n-1}C_{\mathrm{comp}}(G,n-1)
		\sum_{a=1}^d
		\|\partial_a v\|_{H^{n-1}(V)}.
	\end{aligned}
	\]
	In particular, every first weak derivative of \(v\circ G\)
	belongs to \(H^{n-1}(U)\). Together with
	\eqref{eq:l2-composition}, this shows that
	\(v\circ G\in H^n(U)\).
	
	By the Cauchy--Schwarz inequality and the definition of the
	Sobolev norm,
	\[
	\left(
	\sum_{a=1}^d
	\|\partial_a v\|_{H^{n-1}(V)}
	\right)^2
	\le
	d^2\|v\|_{H^n(V)}^2.
	\]
	Furthermore,
	\[
	\|v\circ G\|_{H^n(U)}^2
	\le
	\|v\circ G\|_{L^2(U)}^2
	+
	\sum_{i=1}^d
	\|\partial_i(v\circ G)\|_{H^{n-1}(U)}^2.
	\]
	Combining the preceding estimates with
	\eqref{eq:l2-composition}, we obtain
	\[
	\|v\circ G\|_{H^n(U)}
	\le
	C_{\mathrm{comp}}(G,n)\|v\|_{H^n(V)},
	\]
	where
\[
C_{\mathrm{comp}}(G,n):=
\left(J_G^2+d^3C^2_{\mathrm{mult}}(d,n-1)A_{n-1}^2
[C_{\mathrm{comp}}(G,n-1)]^2\right)^{1/2}.
\]
	Starting with $C_{\mathrm{comp}}(G,1):=J_G\max\{1,L_1\}$ and \eqref{eq:composition-proof-h1},
	induction proves the estimate for every integer \(1\le n\le r\).
\end{proof}

For later use, set $C_{\mathrm{comp}}(G,0):=J_G$.
For $v\in L^2(\Omega)$, define the physical quasi-interpolant by
\begin{equation}
	\label{eq:physical-operator}
	\mathcal Q_{F,h}^\ell v
	:=
	\bigl(\mathcal Q_h^\ell(v\circ F)\bigr)\circ F^{-1}.
\end{equation}
The operator is well defined: \eqref{eq:l2-composition} implies \(v\circ F\in L^2(D)\), and composition with \(F^{-1}\) gives an element of \(L^2(\Omega)\). 

\begin{theorem}[Physical-domain approximation]
	\label{thm:physical-error}
	Assume the hypotheses of Theorem~\ref{thm:error}, with the same
	fixed orders \(\ell,s,m\), and suppose that the fixed geometry map
	\(F:D\to\Omega\) given in \eqref{eq:ext_geo} is bi-Lipschitz.
	If $m\ge1$, assume additionally
	$F^{-1}\in W^{m,\infty}(\Omega;\mathbb{R}^d)$; for $m=0$ bi-Lipschitz
	continuity suffices.
	Then for every $v$ such that $v\circ F\in H^s(D)$,
	\begin{equation}
		\label{eq:physical-pullback-error}
		\|v-\mathcal Q_{F,h}^\ell v\|_{H^m(\Omega)}
		\le
		C_{\mathrm{phys}}^{\mathrm{pb}}\,
		h^{s-m}|v\circ F|_{H^s(D)},
	\end{equation}
	where $C_{\mathrm{phys}}^{\mathrm{pb}}:=C_{\mathrm{comp}}(F^{-1},m)\,C_{\mathrm{par}}$.
	If, in addition, $F\in W^{s,\infty}(D;\mathbb R^d)$,
	then, for every \(v\in H^s(\Omega)\),
	\begin{equation}
		\label{eq:physical-standard-error}
		\|v-\mathcal Q_{F,h}^\ell v\|_{H^m(\Omega)}
		\le
		C_{\mathrm{phys}}^{\mathrm{std}}\,
		h^{s-m}\|v\|_{H^s(\Omega)},
	\end{equation}
	where $C_{\mathrm{phys}}^{\mathrm{std}}:=C_{\mathrm{comp}}(F^{-1},m)\,C_{\mathrm{par}}\,C_{\mathrm{comp}}(F,s)$.
\end{theorem}
\begin{proof}
	Let $u=v\circ F$.
	By hypothesis, \(u\in H^s(D)\). Theorem~\ref{thm:error} gives
	\(u-\mathcal Q_h^\ell u\in H^m(D)\) and
	\[
	\|u-\mathcal Q_h^\ell u\|_{H^m(D)}
	\le
	C_{\mathrm{par}}h^{s-m}|u|_{H^s(D)}.
	\]
	Using \eqref{eq:l2-composition} and Lemma~\ref{lem:composition} with \(G=F^{-1}\) gives
	\[
	\|v-\mathcal Q_{F,h}^\ell v\|_{H^m(\Omega)}=\|\left(u-\mathcal Q_h^\ell u\right)\circ F^{-1}\|_{H^m(\Omega)}
	\le
	C_{\mathrm{comp}}(F^{-1},m)\|u-\mathcal Q_h^\ell u\|_{H^m(D)}.
	\]
	Combining the above two estimates proves \eqref{eq:physical-pullback-error}.
	
	Now suppose \(F\in W^{s,\infty}(D;\mathbb R^d)\) and
	\(v\in H^s(\Omega)\). Lemma~\ref{lem:composition}, applied with
	\(G=F\), gives
	\[
	|v\circ F|_{H^s(D)}
	\le
	\|v\circ F\|_{H^s(D)}
	\le
	C_{\mathrm{comp}}(F,s)\|v\|_{H^s(\Omega)}.
	\]
	Substitution into \eqref{eq:physical-pullback-error} proves
	\eqref{eq:physical-standard-error}.

	Finally, \(\mathcal Q_{F,h}^\ell v\in V_h^\ell\). Moreover,
	\(V_h^\ell\subseteq H^m(\Omega)\) under the stated assumptions.
	Indeed, for \(m\ge1\), this follows from
	Lemmas~\ref{lem:product-space-regularity} and~\ref{lem:composition};
	for \(m=0\), it follows from the \(L^2\) composition estimate.
	Consequently, $\inf_{v_h\in V_h^\ell}\|v-v_h\|_{H^m(\Omega)}\le\|v-\mathcal Q_{F,h}^\ell v\|_{H^m(\Omega)}$.
	The best-approximation estimate follows from the preceding bound.
\end{proof}

\begin{remark}[Space reproduction and operator exactness]
The inclusion $\PP_\ell(\Omega)\subseteq V_h^\ell$ does not imply that
$\mathcal Q_{F,h}^\ell$ reproduces every physical polynomial of degree
at most $\ell$. The operator is exact whenever the pullback belongs to
$\PP_{\ell k}(D)$; it has not been shown to be a projector onto $W_h^\ell$.
For a nonpolynomial geometry, a physical polynomial may have a
nonpolynomial pullback. Moreover, the full norm in
\eqref{eq:physical-standard-error} accounts for lower-order derivatives
in a nonaffine composition.
\end{remark}

\section{Approximation as the product order increases}
\label{sec:pconv}

We retain the product spaces defined in \eqref{eq:intro-product}.
Throughout this section, \(k\ge1\) is fixed and
\begin{equation}\label{pconv:initial-containment}
	\mathbb P_k(D)\subseteq W_h^1.
\end{equation}
This inclusion follows from the polynomial exactness of the
initial operator \(Q_h\).
For each \(\ell\ge1\), set
\[
p_\ell:=\ell k.
\]
The product construction gives
\begin{equation}\label{pconv:product-containment}
	\mathbb P_{p_\ell}(D)\subseteq W_h^\ell.
\end{equation}
Indeed, a monomial of total degree at most \(\ell k\) can be
written as a product of \(\ell\) monomials, each of total degree
at most \(k\), with constant factors included when necessary.

For an integer \(m\ge0\), define
\begin{equation}\label{pconv:best-error}
	E_{h,\ell}^{m}(u)
	:=
	\inf_{w_h\in W_h^\ell\cap H^m(D)}
	\|u-w_h\|_{H^m(D)},
	\qquad u\in H^m(D),
\end{equation}
where \(H^0(D)=L^2(D)\).
If \(W_h^1\subset W^{M,\infty}(D)\) and \(m\le M\),
the weak Leibniz formula gives
\(W_h^\ell\subset H^m(D)\), and the intersection in
\eqref{pconv:best-error} can be omitted.

The parameter \(p_\ell\) measures the guaranteed degree of
complete polynomial reproduction.
At a fixed mesh level, increasing \(\ell\) therefore yields
a sequence of spaces with increasing polynomial
approximation capacity.

\subsection{Polynomial approximation on the parameter domain}
\label{subsec:pconv-polynomials}

We first establish a polynomial approximation estimate with
a constant independent of the polynomial degree.

\begin{lemma}[Complete polynomial approximation]
	\label{pconv:polynomial-lemma}
	Let \(s,m\) be integers satisfying \(0\le m\le s\).
	There is a constant \(C_{d,s,m}>0\) such that, for every
	integer \(p\ge1\) and every \(u\in H^s(D)\),
	\begin{equation}\label{pconv:polynomial-estimate}
		\inf_{\pi\in\mathbb P_p(D)}
		\|u-\pi\|_{H^m(D)}
		\le
		C_{d,s,m}p^{-(s-m)}\|u\|_{H^s(D)}.
	\end{equation}
	The constant is independent of \(p\) and \(u\).
\end{lemma}

\begin{proof}
	The Sobolev extension theorem on the cube provides a bounded
	linear operator
	\[
	\mathcal E_s:H^s(D)\longrightarrow H^s(\mathbb R^d)
	\]
	such that
	\[
	(\mathcal E_su)|_D=u,
	\qquad
	\|\mathcal E_su\|_{H^s(\mathbb R^d)}
	\le C_{\mathrm{ext},d,s}\|u\|_{H^s(D)}.
	\]
	Set \(U=\mathcal E_su\), and use the unitary Fourier transform
	on \(\mathbb R^d\).
	
	For \(R\ge1\), define
	\[
	U_R(\xi)
	=
	(2\pi)^{-d/2}
	\int_{|\eta|\le R}
	\widehat U(\eta)e^{i\xi\cdot\eta}\,d\eta.
	\]
	The integral is well defined because the restriction of
	\(\widehat U\) to a bounded set belongs to \(L^1\).
	By Plancherel's theorem and the equivalence of the
	integer-order Sobolev norm with its Fourier representation,
	\begin{align}
		\|U-U_R\|_{H^m(\mathbb R^d)}^2
		&\le
		C_{d,m}
		\int_{|\eta|>R}
		(1+|\eta|^2)^m|\widehat U(\eta)|^2\,d\eta
		\nonumber\\
		&\le
		C_{d,s,m}R^{-2(s-m)}
		\|U\|_{H^s(\mathbb R^d)}^2.
		\label{pconv:fourier-tail}
	\end{align}
	
	Let
	\[
	\xi_*=(1/2,\ldots,1/2),
	\qquad r_d=\sqrt d/2.
	\]
	For \(p\ge m\), introduce the polynomial
	\begin{equation}\label{pconv:fourier-polynomial}
		P_{p,R}U(\xi)
		=
		(2\pi)^{-d/2}
		\int_{|\eta|\le R}
		\widehat U(\eta)e^{i\xi_*\cdot\eta}
		\sum_{n=0}^{p}
		\frac{\bigl(i(\xi-\xi_*)\cdot\eta\bigr)^n}{n!}
		\,d\eta.
	\end{equation}
	Each summand has total degree at most \(n\), so
	\(P_{p,R}U\in\mathbb P_p\), allowing complex coefficients
	temporarily.
	
	For a multi-index \(\alpha\) with \(|\alpha|\le m\),
	differentiation of the truncated exponential gives
	\[
	\partial^\alpha
	\sum_{n=0}^{p}
	\frac{\bigl(i(\xi-\xi_*)\cdot\eta\bigr)^n}{n!}
	=
	(i\eta)^\alpha
	\sum_{n=0}^{p-|\alpha|}
	\frac{\bigl(i(\xi-\xi_*)\cdot\eta\bigr)^n}{n!}.
	\]
	Since \(|\xi-\xi_*|\le r_d\) on \(D\), the exponential
	series remainder satisfies
	\begin{align}
		&\left|
		\partial^\alpha
		\left(
		e^{i(\xi-\xi_*)\cdot\eta}
		-
		\sum_{n=0}^{p}
		\frac{\bigl(i(\xi-\xi_*)\cdot\eta\bigr)^n}{n!}
		\right)
		\right|
		\nonumber\\
		&\hspace{2cm}\le
		R^{|\alpha|}e^{r_dR}
		\frac{(r_dR)^{p-|\alpha|+1}}
		{(p-|\alpha|+1)!},
		\qquad |\eta|\le R.
		\label{pconv:taylor-remainder}
	\end{align}
	
	Choose \(a_d>0\) sufficiently small that
	\[
	2er_da_d\,e^{2r_da_d}\le\frac12,
	\]
	and set \(R=a_dp\).
	For \(p\ge2m+2\), let \(N=p-|\alpha|+1\).
	Then \(N\ge p/2\). The inequality \(N!\ge(N/e)^N\)
	implies
	\[
	e^{r_da_dp}\frac{(r_da_dp)^N}{N!}
	\le
	\left(2er_da_d\,e^{2r_da_d}\right)^N
	\le 2^{-p/2}.
	\]
	Consequently, for all sufficiently large \(p\), the
	right-hand side of \eqref{pconv:taylor-remainder} is bounded by
	\[
	C_{d,m}p^m2^{-p/2}.
	\]
	
	By Cauchy--Schwarz and Plancherel,
	\[
	\int_{|\eta|\le R}|\widehat U(\eta)|\,d\eta
	\le
	C_dR^{d/2}\|U\|_{L^2(\mathbb R^d)}.
	\]
	Applying this bound to \eqref{pconv:fourier-polynomial},
	and summing over \(|\alpha|\le m\), yields
	\begin{equation}\label{pconv:low-frequency-error}
		\|U_R-P_{p,R}U\|_{H^m(D)}
		\le
		C_{d,m}p^{m+d/2}2^{-p/2}
		\|U\|_{L^2(\mathbb R^d)}.
	\end{equation}
	For fixed \(d,s,m\),
	\[
	\sup_{p\ge1}
	p^{s+d/2}2^{-p/2}<\infty.
	\]
	Combining \eqref{pconv:fourier-tail} and
	\eqref{pconv:low-frequency-error}, with \(R=a_dp\), gives
	\[
	\|u-P_{p,R}U\|_{H^m(D)}
	\le
	C_{d,s,m}p^{-(s-m)}
	\|U\|_{H^s(\mathbb R^d)}.
	\]
	The extension bound proves
	\eqref{pconv:polynomial-estimate} for all sufficiently large
	\(p\). The remaining finitely many degrees are covered by
	taking the zero polynomial and increasing the constant.
	For real-valued \(u\), taking the real part of the
	approximating polynomial preserves the degree and does
	not increase the error.
\end{proof}

\subsection{Algebraic convergence in the product order}
\label{subsec:pconv-algebraic}

\begin{theorem}[Algebraic \(p\)-approximation]
	\label{pconv:algebraic}
	Assume \eqref{pconv:initial-containment}.
	Let \(s,m\) be fixed integers with \(0\le m<s\).
	Then
	\begin{equation}\label{pconv:algebraic-full}
		E_{h,\ell}^{m}(u)
		\le
		C_{d,s,m}p_\ell^{-(s-m)}
		\|u\|_{H^s(D)},
		\qquad
		u\in H^s(D),\quad \ell\ge1.
	\end{equation}
	If \(p_\ell\ge s-1\), then
	\begin{equation}\label{pconv:algebraic-seminorm}
		E_{h,\ell}^{m}(u)
		\le
		C_{d,s,m}p_\ell^{-(s-m)}
		|u|_{H^s(D)}.
	\end{equation}
	The constants are independent of \(h,\ell,k\), and \(u\).
\end{theorem}

\begin{proof}
	By \eqref{pconv:product-containment},
	\[
	\mathbb P_{p_\ell}(D)
	\subseteq W_h^\ell\cap H^m(D).
	\]
	Therefore,
	\[
	E_{h,\ell}^{m}(u)
	\le
	\inf_{\pi\in\mathbb P_{p_\ell}(D)}
	\|u-\pi\|_{H^m(D)}.
	\]
	Lemma~\ref{pconv:polynomial-lemma}, with \(p=p_\ell\),
	proves \eqref{pconv:algebraic-full}.
	
	For the seminorm estimate, we construct
	\(r_u\in\mathbb P_{s-1}(D)\) satisfying
	\begin{equation}\label{pconv:mean-conditions}
		\int_D\partial^\alpha(u-r_u)\,d\xi=0,
		\qquad |\alpha|\le s-1.
	\end{equation}
	The number of these conditions equals
	\(\dim\mathbb P_{s-1}\).
	The corresponding map on \(\mathbb P_{s-1}\) is injective:
	if a polynomial \(r\) satisfies the homogeneous conditions,
	then each derivative of order \(s-1\) is constant and has
	zero integral, hence vanishes. Thus \(r\) has degree at most
	\(s-2\). Repeating the argument gives \(r=0\).
	The map is consequently invertible, proving the existence
	and uniqueness of \(r_u\).
	
	Apply the mean-zero Poincare inequality to each
	\(\partial^\alpha(u-r_u)\), starting with
	\(|\alpha|=s-1\) and proceeding to lower orders.
	Because all derivatives of \(r_u\) of order \(s\) vanish,
	this gives
	\begin{equation}\label{pconv:quotient-estimate}
		\|u-r_u\|_{H^s(D)}
		\le C_{d,s}|u|_{H^s(D)}.
	\end{equation}
	If \(p_\ell\ge s-1\), then
	\(r_u\in\mathbb P_{p_\ell}(D)\). Hence
	\[
	\inf_{\pi\in\mathbb P_{p_\ell}(D)}
	\|u-\pi\|_{H^m(D)}
	=
	\inf_{\pi\in\mathbb P_{p_\ell}(D)}
	\|(u-r_u)-\pi\|_{H^m(D)}.
	\]
	Applying Lemma~\ref{pconv:polynomial-lemma} to \(u-r_u\)
	and using \eqref{pconv:quotient-estimate} proves
	\eqref{pconv:algebraic-seminorm}.
\end{proof}

\begin{corollary}[Density at a fixed mesh level]
	\label{pconv:density}
	Assume \eqref{pconv:initial-containment}.
	For every fixed \(h\), every integer \(m\ge0\), and every
	\(u\in H^m(D)\),
	\[
	\lim_{\ell\to\infty}E_{h,\ell}^{m}(u)=0.
	\]
\end{corollary}

\begin{proof}
	Restrictions to \(D\) of smooth functions on
	\(\mathbb R^d\) are dense in \(H^m(D)\), as follows by
	extending and mollifying.
	Given \(\varepsilon>0\), choose
	\(v\in H^{m+1}(D)\) such that
	\[
	\|u-v\|_{H^m(D)}<\varepsilon.
	\]
	Theorem~\ref{pconv:algebraic}, applied with \(s=m+1\),
	gives
	\[
	E_{h,\ell}^{m}(u)
	\le
	\varepsilon+
	C_{d,m}p_\ell^{-1}\|v\|_{H^{m+1}(D)}.
	\]
	First let \(\ell\to\infty\), and then
	\(\varepsilon\to0\).
\end{proof}

\subsection{Exponential approximation of analytic functions}
\label{subsec:pconv-exponential}

For \(\rho>1\), let \(\mathcal E_\rho\subset\mathbb C\)
be the open region bounded by the Bernstein ellipse
\[
\left\{
\frac12(w+w^{-1}): |w|=\rho
\right\}.
\]
Define the affine map
\[
\Psi(z)=\frac12(z+\boldsymbol{1}),
\qquad z\in\mathbb C^d,
\]
which maps \((-1,1)^d\) onto \(D\).
For a function \(u\) on \(D\), write
\[
\widetilde u=u\circ\Psi.
\]

\begin{lemma}[Analytic polynomial approximation]
	\label{pconv:analytic-lemma}
	Let \(\rho>1\), and suppose that \(\widetilde u\) has a
	holomorphic extension to an open neighborhood of
	\(\overline{\mathcal E_\rho}^{\,d}\).
	Set
	\[
	M_\rho(u)
	:=
	\sup_{z\in\overline{\mathcal E_\rho}^{\,d}}
	|\widetilde u(z)|.
	\]
	For every integer \(m\ge0\) and every
	\(1<\sigma<\rho\), there is a constant
	\(C_{d,m,\rho,\sigma}>0\) such that
	\begin{equation}\label{pconv:analytic-polynomial}
		\inf_{\pi\in\mathbb P_p(D)}
		\|u-\pi\|_{H^m(D)}
		\le
		C_{d,m,\rho,\sigma}M_\rho(u)\sigma^{-p},
		\qquad p\ge1.
	\end{equation}
\end{lemma}

\begin{proof}
	Let \(T_n\) denote the Chebyshev polynomial characterized by
	\[
	T_n(\cos\theta)=\cos(n\theta).
	\]
	The tensor-product Chebyshev expansion of \(\widetilde u\)
	has the form
	\begin{equation}\label{pconv:chebyshev-series}
		\widetilde u(z)
		=
		\sum_{\nu\in\mathbb N_0^d}
		a_\nu\prod_{j=1}^dT_{\nu_j}(z_j).
	\end{equation}
	Its coefficients satisfy
	\begin{equation}\label{pconv:chebyshev-coefficients}
		|a_\nu|
		\le 2^dM_\rho(u)\rho^{-|\nu|}.
	\end{equation}
	
	To obtain this bound, consider
	\[
	g(w)
	=
	\widetilde u\left(
	\frac{w_1+w_1^{-1}}2,\ldots,
	\frac{w_d+w_d^{-1}}2
	\right).
	\]
	This function is holomorphic on a neighborhood of the
	product annulus
	\[
	\rho^{-1}\le|w_j|\le\rho,
	\qquad 1\le j\le d.
	\]
	For nonnegative indices $\nu$, the multivariable Cauchy integral formula bounds the Laurent coefficients by
\[M_\rho(u)\rho^{-|\nu|}.\]
	The identities
	\[
	T_n\left(\frac{w+w^{-1}}2\right)
	=
	\frac{w^n+w^{-n}}2
	\]
	and the invariance of \(g\) under
	\(w_j\mapsto w_j^{-1}\) identify the coefficient \(a_\nu\)
	as \(2^{r(\nu)}\) times the corresponding Laurent
	coefficient, where \(r(\nu)\le d\) is the number of
	nonzero components of \(\nu\).
	This proves \eqref{pconv:chebyshev-coefficients} and
	the expansion \eqref{pconv:chebyshev-series}.
	
	For \(p\ge1\), take the total-degree truncation
	\[
	\widetilde\pi_p(z)
	=
	\sum_{|\nu|\le p}
	a_\nu\prod_{j=1}^dT_{\nu_j}(z_j),
	\qquad
	\pi_p=\widetilde\pi_p\circ\Psi^{-1}.
	\]
	Each summand has total degree \(|\nu|\), so
	\(\pi_p\in\mathbb P_p(D)\).
	
	The univariate Markov inequality, applied repeatedly,
	gives
	\[
	\|T_n^{(r)}\|_{L^\infty(-1,1)}
	\le (n+1)^{2r},
	\qquad r\ge0.
	\]
	For \(r>n\), the left-hand side is zero.
	Consequently, for every multi-index \(\alpha\) with
	\(|\alpha|\le m\),
	\[
	\left\|
	\partial^\alpha
	\prod_{j=1}^dT_{\nu_j}(2\xi_j-1)
	\right\|_{L^\infty(D)}
	\le
	2^{|\alpha|}(1+|\nu|)^{2|\alpha|}.
	\]
	Together with \eqref{pconv:chebyshev-coefficients},
	this also proves absolute and uniform convergence of
	the differentiated series through order \(m\).
	
	There are
	\(\binom{n+d-1}{d-1}\) multi-indices \(\nu\) with
	\(|\nu|=n\). Since \(|D|=1\), it follows that
	\begin{align}
		\|u-\pi_p\|_{H^m(D)}
		&\le
		C_{d,m}M_\rho(u)
		\sum_{n=p+1}^{\infty}
		\binom{n+d-1}{d-1}(n+1)^{2m}\rho^{-n}
		\nonumber\\
		&\le
		C_{d,m}M_\rho(u)
		\sum_{n=p+1}^{\infty}
		(n+1)^{2m+d-1}\rho^{-n}.
		\label{pconv:analytic-tail}
	\end{align}
	For \(1<\sigma<\rho\),
	\[
	C_*:=
	\sup_{n\ge0}
	(n+1)^{2m+d-1}
	\left(\frac{\sigma}{\rho}\right)^n
	<\infty.
	\]
	Therefore,
	\[
	\sum_{n=p+1}^{\infty}
	(n+1)^{2m+d-1}\rho^{-n}
	\le
	C_*\sum_{n=p+1}^{\infty}\sigma^{-n}
	=
	\frac{C_*}{\sigma-1}\sigma^{-p}.
	\]
	Substitution into \eqref{pconv:analytic-tail} proves
	\eqref{pconv:analytic-polynomial}.
	For real-valued \(u\), the real part of \(\pi_p\) gives
	a real polynomial with the same error bound.
\end{proof}

\begin{theorem}[Exponential approximation in the product order]
	\label{pconv:exponential}
	Assume \eqref{pconv:initial-containment} and the analytic
	extension hypothesis of Lemma~\ref{pconv:analytic-lemma}.
	For every fixed integer \(m\ge0\) and every
	\(1<\sigma<\rho\),
	\begin{equation}\label{pconv:exponential-estimate}
		E_{h,\ell}^{m}(u)
		\le
		C_{d,m,\rho,\sigma}M_\rho(u)
		\sigma^{-\ell k}
		=
		C_{d,m,\rho,\sigma}M_\rho(u)
		e^{-b_\sigma\ell k},
		\qquad b_\sigma:=\log\sigma>0.
	\end{equation}
	The constant is independent of \(h,\ell,k\), and \(u\).
\end{theorem}

\begin{proof}
	The inclusion \eqref{pconv:product-containment} gives
	\[
	E_{h,\ell}^{m}(u)
	\le
	\inf_{\pi\in\mathbb P_{\ell k}(D)}
	\|u-\pi\|_{H^m(D)}.
	\]
	Apply Lemma~\ref{pconv:analytic-lemma} with \(p=\ell k\).
\end{proof}

The analytic extension hypothesis holds for some
\(\rho>1\) whenever \(u\) admits a holomorphic extension
to an open complex neighborhood of \(\overline D\).
Indeed, the compact sets
\(\Psi(\overline{\mathcal E_\rho}^{\,d})\) approach
\(\overline D\) as \(\rho\downarrow1\).
Compactness therefore permits a choice of \(\rho>1\)
for which the entire product ellipse lies inside that
neighborhood.

\subsection{Approximation under the fixed geometry map}
\label{subsec:pconv-physical}

Assume that the fixed geometry map in \eqref{eq:ext_geo} is bi-Lipschitz, and retain the physical spaces in \eqref{eq:intro-product}.
For \(u\in H^m(\Omega)\), define
\[
E_{F,h,\ell}^{m}(u)
=
\inf_{v_h\in V_h^\ell\cap H^m(\Omega)}
\|u-v_h\|_{H^m(\Omega)}.
\]

\begin{corollary}[Physical-domain approximation]
	\label{pconv:physical}
	Assume \eqref{pconv:initial-containment}.
	Let \(F:D\to\Omega\) be bi-Lipschitz.
	For \(m\ge1\), assume additionally that
	\[
	F^{-1}\in W^{m,\infty}(\Omega;\mathbb R^d).
	\]
	Then the following statements hold.
	
	\begin{enumerate}
		\item
		If \(s>m\) is an integer and \(u\circ F\in H^s(D)\), then
		\begin{equation}\label{pconv:physical-algebraic}
			E_{F,h,\ell}^{m}(u)
			\le
			C_{F,m}C_{d,s,m}
			(\ell k)^{-(s-m)}
			\|u\circ F\|_{H^s(D)}.
		\end{equation}
		If, in addition,
		\(F\in W^{s,\infty}(D;\mathbb R^d)\), then
		\begin{equation}\label{pconv:physical-algebraic-full}
			E_{F,h,\ell}^{m}(u)
			\le
			C_{F,d,s,m}
			(\ell k)^{-(s-m)}
			\|u\|_{H^s(\Omega)},
			\qquad u\in H^s(\Omega).
		\end{equation}
		
		\item
		If \(u\circ F\) satisfies the analytic extension
		hypothesis of Lemma~\ref{pconv:analytic-lemma}, then
		for every \(1<\sigma<\rho\),
		\begin{equation}\label{pconv:physical-exponential}
			E_{F,h,\ell}^{m}(u)
			\le
			C_{F,m}C_{d,m,\rho,\sigma}
			M_\rho(u\circ F)\sigma^{-\ell k}.
		\end{equation}
	\end{enumerate}
	
	Here \(C_{F,m}\) depends only on \(d,m\), the
	bi-Lipschitz constants, and, for \(m\ge1\), the indicated
	derivative bounds of \(F^{-1}\).
	The constant \(C_{F,d,s,m}\) also depends on the
	\(W^{s,\infty}\) bounds of \(F\).
	All constants are independent of \(h\) and \(\ell\).
\end{corollary}

\begin{proof}
	Let \(\pi\in\mathbb P_{\ell k}(D)\).
	By \eqref{pconv:product-containment},
	\[
	v_h:=\pi\circ F^{-1}\in V_h^\ell.
	\]
	The change-of-variables estimate for \(m=0\), and the
	integer-order Sobolev composition estimate for \(m\ge1\),
	give
	\[
	\|u-v_h\|_{H^m(\Omega)}
	=
	\|((u\circ F)-\pi)\circ F^{-1}\|_{H^m(\Omega)}
	\le
	C_{F,m}
	\|(u\circ F)-\pi\|_{H^m(D)}.
	\]
	In particular, the same composition estimate shows that
	\(v_h\in H^m(\Omega)\).
	Taking the infimum over \(\pi\) yields
	\[
	E_{F,h,\ell}^{m}(u)
	\le
	C_{F,m}
	\inf_{\pi\in\mathbb P_{\ell k}(D)}
	\|(u\circ F)-\pi\|_{H^m(D)}.
	\]
	Lemma~\ref{pconv:polynomial-lemma} proves
	\eqref{pconv:physical-algebraic}.
	When \(F\in W^{s,\infty}(D;\mathbb R^d)\), the
	composition estimate also gives
	\[
	\|u\circ F\|_{H^s(D)}
	\le C_{F,s}\|u\|_{H^s(\Omega)},
	\]
	which proves \eqref{pconv:physical-algebraic-full}.
	Finally, Lemma~\ref{pconv:analytic-lemma}, applied to
	\(u\circ F\), proves \eqref{pconv:physical-exponential}.
\end{proof}

\section{Stable computation in product spaces}
\label{sec:computation}
Fix an integer $\ell\ge1$ and use the spaces in
\eqref{eq:intro-product} and the generators in \eqref{eq:generators}.
Assume \eqref{eq:base-regularity} with $M\ge1$ and that the geometry
map \eqref{eq:ext_geo} is bi-Lipschitz onto a bounded Lipschitz domain.
Lemmas~\ref{lem:product-space-regularity} and~\ref{lem:composition}
then give $V_h^\ell\subset H^1(\Omega)$. The derivative stability
condition \eqref{eq:base-derivative} is not used in the solver analysis.
The construction below realizes a spectral coarse correction in the
spirit of GenEO~\cite{Spillane}, while retaining the product generators.

\subsection{Compatible Galerkin equations}
\begin{assumption}[Coercive energy form]\label{comp:ass:pde}
The measurable symmetric tensor $\mathsf D$ and coefficient $c$ satisfy
$d_*I\le\mathsf D\le d^*I$ and $0<c_*\le c\le c^*<\infty$ almost
everywhere, with $d_*>0$. For $u,v\in H^1(\Omega)$, let
\begin{equation}\label{comp:eq:energy}
a(u,v)=\int_\Omega\bigl((\nabla u)^T\mathsf D\nabla v+cuv\bigr)\,dx,
\qquad \norma v^2=a(v,v),\qquad \mathcal L\in H^1(\Omega)'.
\end{equation}
\end{assumption}
Thus the variational problem has natural boundary conditions and a unique
solution. All integrals, kernel removals, and inverses below are exact.

Retain the nonzero generators through the index set
$\mathcal J_h^\ell=\{\kappa\in\NN^{N_h}:|\kappa|=\ell,
G_{\kappa,h}^\ell\not\equiv0\}$ and define the synthesis operator
\begin{equation}\label{comp:eq:synthesis}
S:\mathbb R^{\mathcal J_h^\ell}\longrightarrow V_h^\ell,
\qquad Sz=\left(\sum_{\kappa\in\mathcal J_h^\ell}
z_\kappa G_{\kappa,h}^\ell\right)\circ F^{-1}.
\end{equation}
For coordinate vectors $e_\kappa$ of this coefficient space, set
$A_{\kappa\eta}=a(Se_\eta,Se_\kappa)$ and
$b_\kappa=\mathcal L(Se_\kappa)$. Surjectivity of $S$ and coercivity give
\begin{equation}\label{comp:eq:kernel}
z^TAz=\norma{Sz}^2,\qquad \ker A=\ker S,\qquad b\in\ran A.
\end{equation}
Thus $Az=b$ determines a unique Galerkin function $u_h=Sz$,
with coefficients determined modulo $\ker S$; this is the standard
redundant-representation setting of~\cite[Section~5.2, preprint version]{HuWinther2021}.

The partition of unity \eqref{eq:partition_unity} gives an explicit
initial-space injection for every product order:
\begin{equation}\label{comp:eq:embed}
(J\mathbf c)_\kappa=\frac1\ell\sum_{i=1}^{N_h}\kappa_i c_i,
\qquad SJ\mathbf c=\left(\sum_{i=1}^{N_h}c_iB_{i,h}\right)\circ F^{-1}.
\end{equation}
To verify the identity, expand
$(\sum_i c_iB_{i,h})(\sum_i B_{i,h})^{\ell-1}$ using the multinomial
normalization in \eqref{eq:generators}. The coefficient of the
$\kappa$th normalized product is $\ell^{-1}\sum_i\kappa_i c_i$.
The initial-space matrix $A_1=J^TAJ$ can be assembled directly as
\[
(A_1)_{ij}=a(B_{j,h}\circ F^{-1},B_{i,h}\circ F^{-1}).
\]
It is positive definite by the basis property in \eqref{eq:initial-space}
and coercivity.

\subsection{Local spectral correction and its stability}
Use one local problem for each original element $K\in\mathcal T_h$.
Let $I_K=\{\kappa\in\mathcal J_h^\ell:
G_{\kappa,h}^\ell|_K\not\equiv0\}$. Here $\mathbb R^{I_K}$ denotes
the space of real families $y=(y_\kappa)_{\kappa\in I_K}$, with its
Euclidean inner product. For a global coefficient vector $z$, set
$z_K=(z_\kappa)_{\kappa\in I_K}$. Define
\[
S_K:\mathbb R^{I_K}\longrightarrow V_h^\ell,\qquad
S_Ky:=\left(\sum_{\kappa\in I_K}y_\kappa G_{\kappa,h}^\ell\right)
\circ F^{-1}.
\]
The associated local space is
\(
\mathcal V_K:=\ran S_K
=\Span\{G_{\kappa,h}^\ell\circ F^{-1}:\kappa\in I_K\}
\subseteq V_h^\ell.
\)
Let $m_\kappa=\#\{K\in\mathcal T_h:\kappa\in I_K\}$ and
$D_K=\operatorname{diag}_{\kappa\in I_K}(m_\kappa^{-1})$; thus
$(D_Ky)_\kappa=y_\kappa/m_\kappa$. Equal distribution of each coefficient gives
\begin{equation}\label{comp:eq:partition}
Sz=\sum_{K\in\mathcal T_h}
\left(\sum_{\kappa\in I_K}\frac{z_\kappa}{m_\kappa}
G_{\kappa,h}^\ell\right)\circ F^{-1}.
\end{equation}
Let $N_K$ and $H_K$ be the energy Gram matrices of the generators
of $S_K$, with the integral in \eqref{comp:eq:energy} taken over
$F(K)$ and $\Omega$, respectively. In particular, \((H_K)_{\kappa\eta}=A_{\kappa\eta},~\forall~\kappa,\eta\in I_K\).
Summing the element energies gives
\begin{equation}\label{comp:eq:sum}
z^TAz=\sum_{K\in\mathcal T_h}z_K^TN_Kz_K,
\qquad N_K,H_K\succeq0.
\end{equation}
Define the full local support regions by
\(
\Omega_K^\star
=\bigcup_{\kappa\in I_K}
\operatorname{ess\,supp}_{\Omega}
\bigl(G_{\kappa,h}^{\ell}\circ F^{-1}\bigr).
\)
Here essential support is taken relative to $\Omega$.
Choose $\rho\ge1$ such that
\begin{equation}\label{comp:eq:overlap}
\#\{K\in\mathcal T_h:x\in\Omega_K^\star\}\le\rho
\qquad\text{for almost every }x\in\Omega.
\end{equation}
Since $\mathcal T_h$ is finite, the above inequality always holds with $\rho=\#\mathcal T_h$.

Set $T_K=D_KH_KD_K$ and $M_K=N_K+T_K$. For local coefficients $y\in\mathbb R^{I_K}$,
$y^TN_Ky$ measures the energy on $F(K)$, whereas
$y^TT_Ky=\norma{S_KD_Ky}^2$ measures the full-domain energy after
coefficient weighting. Since both matrices are positive semidefinite,
$\ker M_K=\ker N_K\cap\ker T_K$; these common-kernel directions
contribute to neither energy. On the Euclidean orthogonal complement
$(\ker M_K)^\perp$, $M_K$ is positive definite, so we solve
\begin{equation}\label{comp:eq:eigenproblem}
T_Ky=\lambda M_Ky,\qquad 0\ne y\in(\ker M_K)^\perp.
\end{equation}
Every eigenpair satisfies
\[
\lambda=\frac{y^TT_Ky}{y^TN_Ky+y^TT_Ky}\in[0,1].
\]
For a fixed $\tau>0$, select the eigenvectors with
$\lambda>\frac{\tau}{1+\tau}$, equivalently $y^TT_Ky>\tau y^TN_Ky$.
By the spectral theorem for symmetric generalized eigenvalue problems
with a positive-definite right-hand matrix, an $M_K$-orthonormal
eigenbasis exists on $(\ker M_K)^\perp$.
Choose an $M_K$-orthonormal eigenbasis $\{y_{K,j}\}_j$ of $(\ker M_K)^\perp$,
with eigenvalues $\lambda_{K,j}$, such that
$y_{K,i}^TM_Ky_{K,j}=\delta_{ij}$. For $q\in\mathbb R^{I_K}$, define
\[
P_Kq=\sum_{j:\,\lambda_{K,j}>\tau/(1+\tau)}
\bigl(y_{K,j}^TM_Kq\bigr)y_{K,j}.
\]
Thus $P_K$ is the $M_K$-orthogonal projection onto the span of the
selected eigenvectors on $(\ker M_K)^\perp$ and is zero on
$\ker M_K$. 
Define $Q_K=I-P_K$; the complementary component $Q_Kq$ retains the
unselected modes and the common-kernel component, its weighted energy is controlled by Lemma~\ref{comp:lem:local} below.
\begin{lemma}[Local spectral estimate]\label{comp:lem:local}
For $q\in\mathbb R^{I_K}$,
\begin{equation}\label{comp:eq:local-bound}
(Q_Kq)^TT_K(Q_Kq)\le\tau q^TN_Kq.
\end{equation}
\end{lemma}
\begin{proof}
 Together with \eqref{comp:eq:eigenproblem}, the $M_K$-orthonormal eigenbasis $\{y_{K,j}\}_j$ of $( \ker M_K)^\perp$ satisfies
\[
y_{K,i}^TT_Ky_{K,j}=\lambda_{K,j}\delta_{ij},\qquad
y_{K,i}^TN_Ky_{K,j}=(1-\lambda_{K,j})\delta_{ij}.
\]
Write an arbitrary $q\in\mathbb R^{I_K}$ as
\[
q=q_0+\sum_{j}\alpha_jy_{K,j},
\]
where $q_0\in\ker M_K$ and $\alpha_j=y_{K,j}^TM_Kq$. Since $\ker M_K=\ker N_K\cap\ker T_K$, we have $N_Kq_0=T_Kq_0=0$.
 The definition of $Q_K=I-P_K$ therefore gives
\[
q^TN_Kq=\sum_j(1-\lambda_{K,j})\alpha_j^2,\qquad
(Q_Kq)^TT_K(Q_Kq)
=\sum_{j:\,\lambda_{K,j}\le\tau/(1+\tau)}
\lambda_{K,j}\alpha_j^2.
\]
For every index in the latter sum,
$\lambda_{K,j}\le\frac{\tau}{1+\tau}$ is equivalent to
$\lambda_{K,j}\le\tau(1-\lambda_{K,j})$. Hence
\[
(Q_Kq)^TT_K(Q_Kq)
\le\tau\sum_{j:\,\lambda_{K,j}\le\tau/(1+\tau)}
(1-\lambda_{K,j})\alpha_j^2
\le\tau\sum_j(1-\lambda_{K,j})\alpha_j^2
=\tau q^TN_Kq.
\]
This proves \eqref{comp:eq:local-bound}.
\end{proof}

Each selected eigenvector with $\lambda_{K,j}>\frac{\tau}{1+\tau}$ defines the global function $S_KD_Ky_{K,j}$.
Collecting these functions over all elements gives the coarse space
\[
\mathcal V_0=\Span\{S_KD_Ky_{K,j}=
\sum_{\kappa\in I_K}
\frac{(y_{K,j})_\kappa}{m_\kappa}G_{\kappa,h}^{\ell}\circ F^{-1}:K\in\mathcal T_h,
\ \lambda_{K,j}>\frac{\tau}{1+\tau}\}.
\]
For its coefficient representation, form one column of $Z$ for each selected pair $(K,j)$: its $\kappa$th entry is
$\frac{(y_{K,j})_\kappa}{m_\kappa}$ for $\kappa\in I_K$ and zero otherwise.
Then $\mathcal V_0=\ran(SZ)$. These columns need not be linearly independent;
the coarse solve below therefore uses a Gram pseudoinverse.

With $r_K=(r_\kappa)_{\kappa\in I_K}$ and $\dagger$ denoting the
Moore--Penrose inverse, define $B_g$ by
adding the coarse correction to the zero-extended local corrections:
\begin{equation}\label{comp:eq:bg}
(B_gr)_\kappa=
\bigl[Z(Z^TAZ)^\dagger Z^Tr\bigr]_\kappa
+\sum_{K:\,\kappa\in I_K}
\bigl[H_K^\dagger r_K\bigr]_\kappa,
\qquad \kappa\in\mathcal J_h^\ell.
\end{equation}
The first term solves the coarse energy system with matrix $Z^TAZ$
and maps the correction back to global coefficients. Each term in the
sum solves the energy system on $\mathcal V_K$ with matrix $H_K$;
its coefficients are extended by zero outside $I_K$ and accumulated.
An empty coarse space contributes zero.
To incorporate an exact correction on the initial space $V_h^1$,
set $B_1=JA_1^{-1}J^T$ and define the preconditioner
\begin{equation}\label{comp:eq:balanced}
B_*=B_1+(I-B_1A)B_g(I-AB_1).
\end{equation}
For a residual $r$, its application consists of
$t=B_1r$, $y=B_g(r-At)$, and $B_*r=t+y-B_1Ay$.
The first step solves the residual equation on $V_h^1$ exactly.
The correction $y$ is computed from the remaining residual, and
$B_1Ay$ removes its initial-space energy projection. Consequently,
$S(y-B_1Ay)$ is energy-orthogonal to $V_h^1$, so the additional
correction preserves the exact initial-space correction.
Moreover, $(I-AB_1)^T=I-B_1A$; since $B_g$ and $B_1$ are symmetric
positive semidefinite, $B_*$ is also symmetric positive semidefinite.

The pseudoinverses are applied by projecting onto the positive eigenspaces
and solving there;
coefficient restriction and zero extension require no explicit matrices.

\begin{theorem}[Stable realization for arbitrary fixed product order]
\label{comp:thm:stability}
Assume \eqref{eq:base-regularity} with $M\ge1$, a bi-Lipschitz
geometry map \eqref{eq:ext_geo}, and Assumption~\ref{comp:ass:pde}.
Let $C_{\rm sd}=2+(2\rho+1)\tau$.
Every $v\in V_h^\ell$ admits $v=v_0+\sum_Kv_K$, with $v_0\in \mathcal V_0$
and $v_K\in \mathcal V_K$, such that
\begin{equation}\label{comp:eq:sd}
\norma{v_0}^2+\sum_K\norma{v_K}^2\le C_{\rm sd}\norma v^2.
\end{equation}
The nonzero eigenvalues of $B_*^{1/2}AB_*^{1/2}$ belong to
$[C_{\rm sd}^{-1},2\rho]$. Moreover, $SB_*A$ acts as the identity
on represented functions in $V_h^1$: if $Sz\in V_h^1$, then
$SB_*Az=Sz$.
\end{theorem}
\begin{proof}
All sums below range over $K\in\mathcal T_h$.
For any coefficient representative $Sz=v$, set
$v_K=S_KD_KQ_Kz_K$ and $v_0=v-\sum_Kv_K$.
By \eqref{comp:eq:partition},
$v_0=\sum_KS_KD_KP_Kz_K\in V_0$.
By the definition of $T_K$, the local estimate
\eqref{comp:eq:local-bound}, and \eqref{comp:eq:sum},
\[
\sum_K\norma{v_K}^2
=\sum_K(Q_Kz_K)^TT_K(Q_Kz_K)
\le\tau\sum_Kz_K^TN_Kz_K
=\tau\norma v^2.
\]
Pointwise Cauchy--Schwarz and \eqref{comp:eq:overlap} yield
$\norma{\sum_Kv_K}^2\le\rho\sum_K\norma{v_K}^2$.
Since $v_0=v-\sum_Kv_K$, it follows that
\[
\norma{v_0}^2+\sum_K\norma{v_K}^2
\le2\norma v^2+(2\rho+1)\sum_K\norma{v_K}^2
\le\bigl(2+(2\rho+1)\tau\bigr)\norma v^2,
\]
which proves \eqref{comp:eq:sd}.

Let $\Pi_0,\Pi_K$ be the energy projections onto $V_0,\mathcal V_K$, and set
$\mathcal T_g=\Pi_0+\sum_K\Pi_K$. The decomposition implies
\[
\norma v^2
=a(v_0,\Pi_0v)+\sum_K a(v_K,\Pi_Kv)
\le\sqrt{C_{\rm sd}}\norma v
\left(\norma{\Pi_0v}^2+\sum_K\norma{\Pi_Kv}^2\right)^{1/2}.
\]
For arbitrary $w_0\in V_0$ and $w_K\in \mathcal V_K$, overlap gives
$\norma{w_0+\sum_Kw_K}^2\le2\rho
(\norma{w_0}^2+\sum_K\norma{w_K}^2)$.
The lower bound below follows from the preceding Cauchy--Schwarz
estimate. For the upper bound, apply the overlap estimate to the
projections and use
$a(\mathcal T_gv,v)=\norma{\Pi_0v}^2+\sum_K\norma{\Pi_Kv}^2$
to obtain
\[
\norma{\mathcal T_gv}^2\le2\rho\,a(\mathcal T_gv,v),\qquad
a(\mathcal T_gv,v)\le\norma{\mathcal T_gv}\norma v.
\]
Combining these inequalities yields
\[
C_{\rm sd}^{-1}\norma v^2\le a(\mathcal T_gv,v)
\le2\rho\norma v^2.
\]
Let $\Pi_1$ be the energy projection onto $V_h^1$. For the orthogonal
split $v=p+q$, where $p=\Pi_1v$ and $q=(I-\Pi_1)v$, the operator
$\mathcal T_* =\Pi_1+(I-\Pi_1)\mathcal T_g(I-\Pi_1)$,
satisfies
\[
a(\mathcal T_*v,v)=\norma p^2+a(\mathcal T_gq,q),\qquad
\norma v^2=\norma p^2+\norma q^2.
\]
Since $C_{\rm sd}^{-1}\le1\le2\rho$, the same lower and upper bounds
hold for $\mathcal T_*$. If $v\in V_h^1$, then $q=0$ and
$\mathcal T_*v=v$.

To identify these function-space operators with the coefficient
preconditioners, let $L$ be any synthesis map and let
$G_{ij}=a(Le_j,Le_i)$ be its energy Gram matrix. Then
$\ker G=\ker L$. For $g=(a(v,Le_i))_i$, one has
$g\perp\ker G$, and hence $GG^\dagger g=g$.
Consequently, $LG^\dagger g\in\ran L$ satisfies
$a(v-LG^\dagger g,Lc)=0$ for every coefficient vector $c$.
It is therefore the energy projection of $v$ onto $\ran L$,
even when the generators are dependent.
Apply this observation to $SZ$, $S_K$, and $SJ$ to obtain
$SB_gA=\mathcal T_gS$ and $SB_*A=\mathcal T_*S$.
In energy-orthonormal coordinates, write $S$ as $\mathcal S$.
Then $A=\mathcal S^T\mathcal S$, and $\mathcal T_*$ is represented by
$\mathcal S B_*\mathcal S^T$. This matrix and
$B_*^{1/2}\mathcal S^T\mathcal S B_*^{1/2}$ have the same nonzero
eigenvalues. The preceding energy bounds complete the proof.
\end{proof}

\subsection{Algorithm, approximation error, and scope}
For $r=b-Az$, \eqref{comp:eq:kernel} and
Theorem~\ref{comp:thm:stability} imply
\begin{equation}\label{comp:eq:residual}
\norma{u_h-Sz}^2=r^TA^\dagger r
\le C_{\rm sd}\,r^TB_*r.
\end{equation}
Indeed, write $r=Ae$ and apply the lower energy bound to $Se$.
Table~\ref{comp:tab:algorithm} gives compatible preconditioned conjugate
gradients. The spectral bound makes $B_*$ positive definite on
$\ran A$; equivalently, $A^{1/2}B_*A^{1/2}$ is positive definite on
$\ran A$. Thus the recurrences are ordinary positive-definite CG on
the coefficient quotient, even though the stored matrices may be singular.

\begin{table}[htbp]
\centering
\caption{Spectrally stabilized product-space Galerkin solve, for any fixed $\ell\ge1$.}
\label{comp:tab:algorithm}
\small
\begin{tabularx}{\textwidth}{@{}r >{\raggedright\arraybackslash}X@{}}
\toprule
Step & Operation \\
\midrule
Input & Initial basis values, gradients, and element-support incidence;
the data in
\eqref{comp:eq:energy}, product order $\ell$, threshold $\tau>0$,
and algebraic tolerance $\varepsilon>0$. \\
1 & Enumerate $\mathcal J_h^\ell$; use \eqref{eq:generators} and
the product rule to assemble $A,b,N_K,H_K$ on the original mesh.
Form $J$ by \eqref{comp:eq:embed}; assemble $A_1$ directly from the
initial basis. \\
2 & Form $D_K,T_K,M_K$. Restrict $M_K$ to its positive eigenspace,
solve \eqref{comp:eq:eigenproblem}, and retain modes with
$\lambda>\tau/(1+\tau)$. Form $Z$ by zero extension and prepare the local, coarse, and
initial-space solves; choose a verified bound $\rho$ satisfying
\eqref{comp:eq:overlap}. \\
3 & To apply $w=B_*r$, compute
$t=B_1r$, $y=B_g(r-At)$, and $w=t+y-B_1Ay$. \\
4 & Set $z=0$, $r=b$, $w=B_*r$, $\delta=r^Tw$, and $p=w$. \\
5 & While $C_{\rm sd}\delta>\varepsilon^2$, perform steps 6--7. \\
6 & Set $\alpha=\delta/(p^TAp)$;
update $z\leftarrow z+\alpha p$ and $r\leftarrow r-\alpha Ap$. \\
7 & Set $w=B_*r$, $\delta_{\rm new}=r^Tw$,
$p\leftarrow w+(\delta_{\rm new}/\delta)p$, and
$\delta\leftarrow\delta_{\rm new}$. \\
Output & $Sz\in V_h^\ell$, with $\norma{u_h-Sz}\le\varepsilon$. \\
\bottomrule
\end{tabularx}
\end{table}

If Theorem~\ref{thm:physical-error} also applies with $m=1$ and
$1<s\le\ell k+1$, Galerkin best approximation and
\eqref{comp:eq:residual} give
\begin{equation}\label{comp:eq:total-error}
\norma{u-Sz}^2\le
C^2h^{2(s-1)}|u\circ F|_{H^s(D)}^2+\varepsilon^2,
\end{equation}
where $C=\max\{d^*,c^*\}^{1/2}C_{\mathrm{phys}}^{\mathrm{pb}}$.
Here $u$ solves \eqref{comp:eq:energy} with load $\mathcal L$;
the squared splitting follows from Galerkin orthogonality.

The same construction applies to $\ell>2$ without a new spectral proof:
\eqref{eq:generators} supplies the products and \eqref{comp:eq:embed}
supplies their initial-space injection. If at most $q$ initial basis
functions are active per cell and each meets at most $n_{\mathrm{supp}}$
cells, then at most
$q_\ell=\binom{q+\ell-1}{\ell}$ products are active per cell.
Each product meets at most $n_{\mathrm{supp}}$
cells, so one may take $\rho=q_\ell n_{\mathrm{supp}}$.
With fixed $\ell,q,n_{\mathrm{supp}}$, uniformly bounded integration cost,
and $n_0$ selected coarse columns, dense local and coarse factorizations give
\[
\text{setup}=O(N_{\rm cell}+n_0^3)+C_{1,\rm setup},\qquad
\text{application}=O(N_{\rm cell}+n_0^2)+O(C_{1,\rm solve}),
\]
where $N_{\rm cell}=\#\mathcal T_h$ and the last terms are the costs
of the initial-space solver. Neither the overlap nor these costs are
asserted to be uniform in $\ell$, and $n_0$ need not be small.

Bernstein extraction is optional for evaluating the products; it is not
used in Theorem~\ref{comp:thm:stability}. The product rule permits direct
evaluation from the initial basis, including at zeros of a basis function,
without dividing by basis values. For positive-weight quadrature with
positive semidefinite local energy contributions and a coercive assembled
form, the argument can be repeated with the overlap measured at quadrature
nodes. Consistency with \eqref{comp:eq:energy} requires a separate estimate.
Essential boundary conditions require a conforming constrained synthesis
and an injection mapping the constrained initial space into it. Removing
generators or replacing small nonzero eigenvalues by zero is not covered
by the present exact-kernel argument.

\section{Conclusions}
Product enrichment under a fixed geometry map gives simultaneous
parametric and physical polynomial containment. The local construction
turns the initial reproduction and stability assumptions into a
quasi-interpolant with higher polynomial exactness and fixed-order
Sobolev approximation estimates. Derivative and physical-domain estimates
use the regularity conditions stated in their respective theorems.
Increasing-order best approximation follows separately from polynomial
containment.

For every fixed product order, the coefficient synthesis and exact
initial-space injection support a redundant Galerkin formulation with
a two-level energy bound. The resulting solver can be coupled to the
approximation estimate through the preconditioned residual. The cost
bounds retain the coarse and initial solves explicitly. Uniform
multilevel complexity, quadrature consistency, and numerical rank
truncation remain subjects for further analysis.

\end{document}